\documentclass{amsart}
\usepackage{amssymb, amscd, latexsym, graphicx, psfrag}

\newtheorem{dummy}{dummy}[section]

\newtheorem{theorem}[dummy]{Theorem}

\theoremstyle{definition}

\newtheorem{example}[dummy]{Example}
\newtheorem{remark}[dummy]{Remark}

\newcommand{\bC}{\mathbb{C}}

\newcommand{\bF}{\mathbb{F}}
\newcommand{\bL}{\mathbb{L}}
\newcommand{\bP}{\mathbb{P}}
\newcommand{\bQ}{\mathbb{Q}}
\newcommand{\bR}{\mathbb{R}}
\newcommand{\bZ}{\mathbb{Z}}

\newcommand{\cB}{\mathcal{B}}
\newcommand{\cC}{\mathcal{C}}

\newcommand{\cG}{\mathcal{G}}

\newcommand{\cP}{\mathcal{P}}
\newcommand{\cQ}{\mathcal{Q}}
\newcommand{\cO}{\mathcal{O}}

\newcommand{\cT}{\mathcal{T}}

\newcommand{\ft}{\mathfrak{t}}

\newcommand{\Hom}{\mathrm{Hom}}

\newcommand{\si}{\sigma}
\newcommand{\Si}{\Sigma}

\newcommand{\tY}{\tilde{Y}}

\newcommand{\vc}{ {\vec{c}} }

\newcommand{\vt}{ {\vec{t}} }

\newcommand{\tri}{\triangle}
\newcommand{\LS}{ {\Lambda_\Sigma} }

\newcommand{\Spec}{\mathrm{Spec}}
\newcommand{\Ext}{\mathrm{Ext}}

\newcommand{\Sh}{\mathit{Sh}}

\newcommand{\Perf}{\cP\mathrm{erf}}

\begin{document}

\title[Coherent-Constructible Correspondence and HMS for
Toric Varieties]{The Coherent-Constructible Correspondence and
Homological Mirror Symmetry
for Toric Varieties}
\author{Bohan Fang}
\address{Bohan Fang, Department of Mathematics, Northwestern University,
2033 Sheridan Road, Evanston, IL  60208}
\email{b-fang@math.northwestern.edu}

\author{Chiu-Chu Melissa Liu}
\address{Chiu-Chu Melissa Liu, Department of Mathematics, Columbia University,
2990 Broadway, New York, NY 10027}
\email{ccliu@math.columbia.edu}

\author{David Treumann}
\address{David Treumann, School of Mathematics, University of Minnesota,
127 Vincent Hall, 206 Church St. S.E., Minneapolis, MN 55455}
\email{treumann@math.umn.edu}

\author{Eric Zaslow}
\address{Eric Zaslow, Department of Mathematics, Northwestern University,
2033 Sheridan Road, Evanston, IL  60208}
\email{zaslow@math.northwestern.edu}

\dedicatory{Dedicated to Shing-Tung Yau on the occasion of his 59th birthday}

\maketitle

\tableofcontents

\section{Introduction}

Toric geometry has always been closely intertwined with combinatorics.
Toric geometry has also figured prominently in mirror symmetry, particularly
in the derivation of Hori and Vafa, who define the mirrors of toric varieties
and hypersurfaces therein.  (Indeed, toric geometry is responsible for
the vast majority of examples of mirror symmetry.)
In this survey, we give an exposition without proofs of the
results of \cite{NZ, N, FLTZ}, which together link up these
various relations.

Briefly, mirror symmetry relates
coherent sheaves on a toric variety to a Fukaya-type category of Lagrangian submanifolds of
an affine space\footnote{Several possible versions of this Fukaya category
arise in the literature -- see Section \ref{subsec:fs}, \ref{subsec:hirz}
and References \cite{S,Ab2,FSS} --
though we mainly work with the one defined in \cite{NZ}.}
$(\bC^*)^n\cong T^*((S^1)^n).$  An equivariant version of this duality
relates equivariant coherent sheaves to a Fukaya category of
$T^*(\bR^n),$ where $\bR$ is the universal cover of $S^1.$
The thereoms of \cite{NZ,N} equate
the Fukaya category of a cotangent $T^*B$ to the category of constructible sheaves
on the base $B.$  In our application $B\cong \bR^n$ is a vector space and the
constructible category
is generated by constant sheaves on polytopes.   It is thus a combinatorial
category which completes a triangle of linkages
between equivariant coherent sheaves on a toric variety,
a Fukaya category on a cotangent space of a vector space,
and a constructible sheaf category
generated by polytopes on that vector space.
The resulting map from coherent to constructibe sheaves
is a categorification of Morelli's combinatorial description of
the equivariant K-theory of a toric variety \cite{M}.

\subsection{Outline}

In Section \ref{sec:ms} we discuss the categories arising in
homological mirror symmetry and give a list of some
results to date.  In Section \ref{sec:tduality},
we review the toric geometry of the B-side and the different geometries of
the A-side.  We then describe
how, beginning with an equivariant ample line bundle on a toric
variety $X_\Sigma$ constructed from a fan $\Sigma\subset N_\bR,$
we construct a Lagrangian submanifold of the T-dual mirror
geometry.  This defines a functor from coherent
sheaves on a toric variety
to a Fukaya category of a cotangent bundle.
Section \ref{sec:micro} is a largely self-contained
summary of the equivalence constructed in \cite{NZ,N}, and
can therefore be read (or skipped) independently.  This
section relates the Fukaya category of a cotangent bundle
to constructible sheaves on the base manifold.
In Section \ref{sec:ccc} we review the direct map from
equivariant coherent sheaves on $X_\Sigma$
to constructible sheaves on a vector space $M_\bR = N_\bR^\vee,$
and discuss the relation to Morelli's theorem.
Finally, we try to make these results as accessible as possible
by providing a selection of examples in Section \ref{sec:exs}.

\subsection*{Acknowledgments}
Eric and ML are greatly indebted to
Professor Shing-Tung Yau for his guidance and support over the years.
We thank him for his many kindnesses.
\medskip

The work of EZ is supported in part by NSF/DMS-0707064.
The work of ML is supported in part by Sloan Research
Fellowship.
BF and EZ would like to thank the Pacific Institute for
Mathematical Sciences, where some of this work was performed.

\section{Mirror Symmetry for Toric Manifolds}
\label{sec:ms}

Mathematical treatments of mirror symmetry can take different
forms.  We will focus on the one most suited to our purpose,
the homological mirror symmetry between a
projective toric variety and
its Hori-Vafa Landau-Ginzburg mirror.

\subsection{Hori-Vafa mirror}\label{sec:HV}
In \cite{HV}, Hori and Vafa derived from physics that the mirror geometry
of a projective toric Fano manifold $X$ is a Landau-Ginzburg model
$((\bC^*)^n, W)$, where $n = \dim_\bC X$, and
$W:(\bC^*)^n \to \bC$ is a holomorphic function known as the
superpotential \cite{HV}.
Their method involves T-duality for the real torus acting on $X.$
On $X$ one can consider the A-model $A(X,\omega),$ which only depends on
the symplectic structure $\omega$ of $X$, and the B-model
$B(X,J)$ which only depends on the complex structure $J$ of $X$.
One is then interested in the mathematical structures that arise in
these models. For example, quantum cohomology and
Lagrangian submanifolds arise in the A-model, while
singularity theory and
coherent sheaves arise in the B-model. Schematically,
mirror symmetry postulates the following equivalences:
\begin{eqnarray}
\label{amod} A(X,\omega) \cong B((\bC^*)^n,W),\\
\label{bmod} B(X,J) \cong A((\bC^*)^n, W).
\end{eqnarray}

In \cite{Au1}, Auroux described a construction of the mirror of a non-toric
Fano manifold, or more generally, a K\"{a}hler manifold with a nonzero,
effective anti-canonical divisor.

\subsection{Categories in mirror symmetry}
Kontsevich's homological mirror symmetry conjecture \cite{Ko}
relates categories of D-branes in mirror-dual physical theories.\footnote{The
conjecture as originally posed in \cite{Ko} involved only mirror symmetry of
Calabi-Yau manifolds, but Kontsevich later generalized the conjecture
to the toric Fano case.}
We now introduce the categories
that arise in $A(X,\omega)$, $B(X,J)$,
$A((\bC^*)^n, W)$, and $B((\bC^*)^n,W)$.

\subsubsection{The Fukaya category  $Fuk(X)$}\label{sec:fukaya}
The category that arises in $A(X,\omega)$
is the Fukaya category $Fuk(X)$.  The Fukaya category is a rather vast
subject in general.  We give a superficial treatment
here, and a somewhat expanded discussion in Section \ref{sec:micro} --
see \cite{FOOO,S} for foundational material.

An object in $Fuk(X)$ is a closed
Lagrangian
submanifold of $(X,\omega)$. The hom space of two
Lagrangian submanifolds $L$, $L'$
is the Floer complex:
\begin{equation}\label{eqn:floer-complex}
hom_{Fuk(X)}(L, L') = CF^*(L, L') =\bigoplus_{p\in L\cap \phi(L')}
\Lambda_{nov}\, p
\end{equation}
where $\phi$ is a Hamiltonian symplectic automorphism of
$(X,\omega)$ such that $L$ and
$\phi(L')$ intersect transversally, and
$$
\Lambda_{nov}=\Bigl\{ \sum_{i=0}^\infty a_i T^{\lambda_i} e^{n_i}\mid
a_i\in \bQ, \lambda_i\in \bR_{\geq 0}, n_i\in \bZ,
\lim_{i\to \infty} \lambda_i =\infty\Bigr\}
$$
is the Novikov ring (we will see later that
the sequence $\lambda_i$ consists of possible symplectic areas
of certain holomorphic  polygons in $X$). $Fuk(X)$ is
an $A^\infty$-category:
there are maps
$$
m_k: hom_{Fuk(X)}(L_0, L_1)\otimes \cdots \otimes
hom_{Fuk(X)}(L_{k-1}, L_k) \to hom_{Fuk(X)}(L_0, L_k)[k-2]
$$
satisfying certain bilnear relations.  In particular,
$$
d:=m_1: hom_{Fuk(X)}(L_0,L_1)\to hom_{Fuk(X)}(L_0,L_1)[-1]
$$
is the differential of the Floer complex,\footnote{There may be a ``curved''
obstruction, meaning $d^2\neq 0.$}
and $\circ := m_2:
hom_{Fuk(X)}(L_0,L_1)\otimes hom_{Fuk(X)}(L_1,L_2)\to hom_{Fuk(X)}(L_0,L_2)$
is a (not necessarily associative) compatible product:
$d(p\circ q) = (dp)\circ q \pm p \circ (dq).$

The composition maps $m_k$ are defined by summing over
holomorphic maps $u$ satisfying boundary conditions
as in Figure 1 (assuming $L_i$'s intersect transversally),  weighted by
$$
T^{\textup{ symplectic area of $u$ }} e^{\textup{ Maslov index of $u$ }} \in \Lambda_{nov}.
$$
\begin{figure}[h]
\begin{center}
\psfrag{L0}{\small $L_0$}
\psfrag{L1}{\small $L_1$}
\psfrag{Lk}{\small $L_k$}
\psfrag{p1}{\small $p_1$}
\psfrag{p2}{\small $p_2$}
\psfrag{pk}{\small $p_k$}
\psfrag{q}{\small $q$}
\psfrag{x0}{\small $x_0$}
\psfrag{x1}{\small $x_1$}
\psfrag{x2}{\small $x_2$}
\psfrag{xk}{\small $x_k$}
\psfrag{u}{\small $u$}
\psfrag{ux=q}{\small $u(x_0)=q$}
\psfrag{ux=p}{\small $u(x_i)=p_i,\quad i>0$}
\includegraphics[scale=0.6]{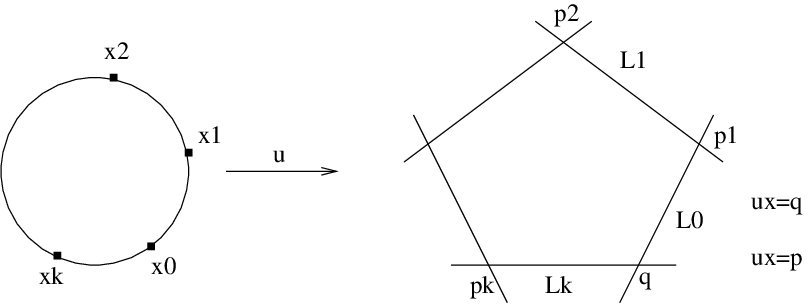}
\end{center}
\caption{Holomorphic polygons}
\end{figure}

More precisely, note that the symplectic area
and Maslov index of $u$ is determined
by its relative homotopy class
$$
\beta \in \pi_2:= \pi_2(X, L_0\cup \cdots\cup L_k).
$$
Let $\int_\beta \omega$ and $\mu(\beta)$
denote the symplectic area
and the Maslov index of the class $\beta$, and
let $n_\beta(p_1,\ldots,p_k,q)$ be the number
of holomorphic maps in Figure 1 in class $\beta$
(sometimes this is a ``virtual''
number which is a rational number instead
of an integer). Then
\begin{equation}\label{eqn:mk}
m_k(p_1,\ldots,p_k) = \sum_q \sum_{\beta\in \pi_2}
n_\beta(p_1,\ldots, p_k,q)
T^{\int_\beta \omega}e^{\mu(\beta)} \cdot q.
\end{equation}
See \cite{FOOO} for details.

\begin{remark}\label{exact-Fuk}
In this paper, we are more interested in {\em exact} Fukaya categories
for a noncompact manifold $X$ equipped  with an exact symplectic form
$\omega = d\theta$, where $\theta$ is a 1-form on $X$.
A Lagrangian submanifold of $(X,\omega)$ is exact if $\theta\vert_L = df$
for some function $f:L\to \bR$. We build the category from exact, possibly
noncompact Lagrangian submanifolds (satisfying certain regularity
properties near infinity). In this case, we my use $\bQ$ or $\bC$
instead of the Novikov ring $\Lambda_{nov}$ in \eqref{eqn:floer-complex},
and define
$$
m_k(p_1,\ldots,p_k) = \sum_q n(p_1,\ldots, p_k,q) \cdot q
$$
where $n(p_1,\ldots,p_k,q)$ is the number of holomorphic maps
satisfying boundary conditions as in Figure 1.
See \cite{S} for details.
\end{remark}

\begin{remark}
In general Lagrangians are equipped with flat complex line bundles.
On the right hand side of \eqref{eqn:floer-complex}, $\Lambda_{nov}$
is replaced by $\Lambda_{nov}\otimes \Hom_{\bC}(V_p, V'_p) \cong
\Lambda_{nov}\otimes_{\bQ}\bC$, where $V$ and $V'$ are flat complex
line bundles on $L$ and $L'$, respectively. The right hand side of
\eqref{eqn:mk} contains an additional factor which is the holonomy
of the flat line bundle along $\partial \beta$. These extra data can
be ignored in the exact case described in the above Remark
\ref{exact-Fuk}.
\end{remark}

The bounded derived category $DFuk(X)$ is
a triangulated category obtained
by taking cohomology $H^0$ of the triangulated
envelope $Tr Fuk(X)$ of $Fuk(X)$ (see Section
\ref{sec:algebra}).

\subsubsection{The dg category of coherent sheaves $Coh(X)$}
The category that arises in $B(X,J)$
is the derived category (or rather an appropriate dg enhancement) of coherent sheaves on $X$.
The shortest way to define $D Coh(X)$ is as the dg category whose objects are complexes of
injective quasicoherent sheaves $I^\bullet$ satisfying
\begin{enumerate}
\item[] The cohomology sheaves of $I^\bullet$ are coherent, and vanish in all but finitely many degrees.
\end{enumerate}
Then $hom(I^\bullet,J^\bullet)$ is the usual chain complex of homomorphisms.

\subsubsection{The Fukaya-Seidel category $FS( (\bC^*)^n, W)$}
\label{subsec:fs}
The category of D-branes that arises in $A((\bC^*)^n,W)$ is the Fukaya-Seidel
category $FS((\bC^*)^n,W)$.

Assume that the superpotential $W$ has isolated
nondegenerate critical points $x_1,\ldots,x_v$,
so that $W:(\bC^*)^n\to \bC$ is a Lefschetz fibration.
Let $b\in \bC$ be a regular value, so that
the fiber $W^{-1}(b)$ is a smooth complex
hypersurface in $(\bC^*)^n$.
An object in $FS((\bC^*),W)$ is a
vanishing cycle which is a Lagrangian
sphere in $W^{-1}(b)$ associated to a path
from $b$ to a critical point.
Alternatively,
it is a Lagrangian thimble
which is a Lagrangian submanifold
of $(\bC^*)^n$ with boundary in
$W^{-1}(b)$; its image under $W$
is a path from a critical value
to $b$.
As in Section \ref{sec:fukaya}, the hom space
of two vanishing cycles
is a Floer complex, and the composition
maps are defined by counting holomorphic
polygons in $W^{-1}(b)$ with boundaries
in vanishing cycles.

We write $DFS((\bC^*)^n,W)$ for the bounded derived
category of the Fukaya-Seidel category
$FS((\bC^*)^n,W)$.

\subsubsection{$DSing((\bC^*)^n, W)$ }
\label{subsec:dsing}
The derived category that arises in the Landau-Ginzburg B-model $B((\bC^*)^n,W)$ is
$D Sing((\bC^*)^n,W)$.

We briefly describe the definition, following
Orlov \cite{Or}.
Assume that the superpotential $W$ has isolated
nondegenerate critical points $x_1,\ldots,x_v$.
Let $W_i$ be the closed subscheme of $(\bC^*)^n$
defined by $W(z)-x_i=0$. Define the
triangulated category $DSing(W_i)$ as the
quotient of the bounded
derived cateogry of coherent sheaves
on $W_i$, $DCoh(W_i)$, by the
full triangulated subcategory
generated by perfect complexes, $\Perf(W_i)$.
Then
$$
DSing((\bC^*)^n,W) =\prod_{i=1}^v DSing(W_i).
$$
Orlov, using
the theorem of Eisenbud (Section 5 of \cite{E}), relates this
category to the category of matrix factorizations.  We will
not use this description or say anything more about it, however.

\subsubsection{Kontsevich's conjecture}

The homological mirror conjecture
for toric Fano manifolds postulates the following
quasi-equivalences of triangulated categories:
\begin{enumerate}
\item[(a)] \label{first} $D Fuk(X)\cong D Sing((\bC^*)^n, W)$.
\item[(b)] \label{second} $D Coh(X) \cong D FS((\bC^*)^n, W)$.
\end{enumerate}
where items (a) and (b) correspond to
Equations (\ref{amod}) and (\ref{bmod}) in Section \ref{sec:HV}, respectively.

\subsection{Results to date}
Here is a run-down of some results in
establishing Kontsevich's conjecture in the
special case of a non-Calabi-Yau toric variety.
Many other resuts in mirror symmetry are omitted, with apologies.

The equivalence (a) has been studied by Cho \cite{Ch1, Ch2},
Cho-Oh \cite{CO}; Fukaya-Oh-Ohta-Ono studied
(a) for any projective toric manifolds, including
non-Fano ones \cite{FOOO1, FOOO2}.

The equivalence (b) has been studied in the following works:
\begin{enumerate}

\item[(i)] Hori, Iqbal, and Vafa define a correspondence
of branes for statement (b) in the
case of projective spaces and toric del Pezzo surfaces in \cite{HIV}.

\item[(ii)] Seidel proves (b) for the complex
projective plane  \cite{S2.5}.

\item[(iii)] In \cite{AKO}, Auroux-Katzarkov-Orlov prove (b) for weighted projective planes and their non-commutative deformations.
They also study Hirzebruch surfaces $\bF_m =\bP(\cO_{\bP^1}\oplus
\cO_{\bP^1}(m))$ ($m\geq 0$), which are not Fano when $m\geq 2$. In
\cite{AKO2}, they consider the blow-ups of $\bP^2$, and prove (b)
for del Pezzo surfaces (including non-toric ones).

\item[(iv)] In \cite{U, UY}, Ueda and Ueda-Yamazaki prove (b) for toric del Pezzo surfaces and their orbifold versions.

\item[(v)] Bondal-Ruan announced a proof of (b)
for weighted projective spaces of any dimension.

\item[(vi)] In \cite{Ab2}, Abouzaid studied the equivalence (b)
for projective toric manifolds, including non-Fano ones.
Abouzaid established a quasi-equivalence between
$D Coh(X)$ and a full subcategory  of $D^\pi Fuk( (\bC^*)^n,W^{-1}(0))$,
where $Fuk( (\bC^*)^n, W^{-1}(0))$ is certain relative
Fukaya category,
and the superscript $\pi$ stands for the split closure.
An object in $Fuk((\bC^*)^n, W^{-1}(0))$ is a
compact Lagrangian submanifold  in $(\bC^*)^n$ with boundary
in the complex hypersurface $W^{-1}(0)$.

\item[(vii)] In \cite{Fa}, the first author proved a version of (b)
for projective spaces $\bP^n$ using T-duality.
On the right
hand side of (b), he uses
a Fukaya category of the contangent bundle
$T^*(S^1)^n \cong (\bC^*)^n$ constructed
by Nadler and the fourth author in \cite{NZ} (see Section \ref{sec:micro}).
The quasi-equivalence is established via
a category of constructible sheaves on
the torus $(S^1)^n$.
\end{enumerate}

We will give an exposition of the recent work by the
authors \cite{FLTZ}  generalizing (v) to all projective
toric manifolds, including non-Fano ones.

\section{T-duality} \label{sec:tduality}

In \cite{SYZ}, mirror symmetry for Calabi-Yau
manifolds was described by dualizing the fibers
of a conjectural special Lagrangian torus fibration with singularities.
This ``T-duality'' provides a transformation of branes.
This transformation
was investigated in \cite{AP,LYZ},
and can be applied to the open orbit $Y\cong (\bC^*)^n$
inside $X$
This is the method we will employ to construct Lagrangian objects
from holomorphic ones.

To be more explicit, we need to introduce
some notation.

\subsection{Moment Polytope}
Here we recall the moment map construction, since it
figures prominently in the sequel.  We ignore many other
constructions of toric geometry.
Readers may consult \cite{Fu} for foundations on toric
geometry, or \cite[Section 2]{FLTZ} for material
immediately related to present purposes.

Let $X$ be a projective toric manifold, and let
$T\cong(\bC^*)^n$ act on $X$, where $n$ is the complex
dimension of $X$. Let $N=\Hom(\bC^*,T)\cong \bZ^n$
be the lattice of 1-parameter subgroups of $T$,
and let $M=\Hom(T,\bC^*)$ be the group
of irreducible characters of $T$. Note that
$M$ is the dual lattice of $N$. Let
$T_\bR \cong U(1)^n$ be the maximal compact
subgroup of $T$. Define
$$
N_\bR = N\otimes _\bZ \bR,\quad
M_\bR = M\otimes _\bZ \bR.
$$
Then $N_\bR$ is canonically identified with
the Lie algebra $\ft_\bR$ of $T_\bR$, while
$M_\bR$ is canonically identified with
$\ft_\bR^*$, the dual real vector space
of $\ft_\bR$. We have
$$
N_\bR/N \cong  T_\bR,\quad
M_\bR/M \cong  T_\bR^\vee
$$
where $T_\bR^\vee$ is the dual torus.

Let $\omega$ be a $T_\bR$-invariant
symplectic form on $X$ such that the $T_\bR$-action
on $X$ is Hamiltonian. Let
$\mu_{T_\bR}: X\to M_\bR$ be the moment map
of the $T_\bR$-action.
It is determined by $\omega$ up to
addition of a constant vector in $M_\bR$.
The image of the moment map is a convex
polytope $\tri\subset M_\bR$, known
as the moment polytope.  The fibers of the
moment maps are $T_\bR$-orbits. We have
a homeomorphism $X/T_\bR \cong \tri$, and
the natural projection $X\to \tri$ restricts
to a $T_\bR$-fibration from the
open orbit $Y\cong (\bC^*)^n$ to the interior
$\tri^\circ$ of the convex polytope $\tri$.

\begin{example}
$X=\bP^2$, and $\omega$ is the Fubini-Study
K\"{a}hler form, so that $\int_{\bP^2} \omega^2 = 4\pi^2$.
The moment map $\mu_{T_\bR}: \bP^2\to M_\bR \cong \bR^2$
is given by
$$
\mu_{T_\bR}([X_0,X_1,X_2])
=\frac{(|X_1|^2, |X_2|^2)}{|X_0|^2 + |X_1|^2 + |X_2|^2}.
$$
The moment polytope is the triangle:
\begin{figure}[h]
\psfrag{(0,0)}{\small $(0,0)$}
\psfrag{(1,0)}{\small $(1,0)$}
\psfrag{(0,1)}{\small $(0,1)$}
\begin{center}
\includegraphics[scale=0.6]{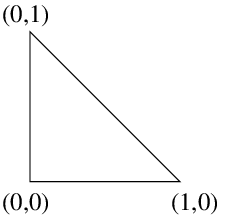}
\end{center}
\caption{Moment polytope of $\bP^2$}
\end{figure}
\end{example}

\subsection{Geometry of the open orbit}
Here we discuss the various geometries arising from the
open orbit of the complex torus under T-duality.\footnote{We use the convention
that the dual of a circle of radius $R$ is a circle of radius
$1/R$. An alternative convention is that the dual
of a circle of {\em circumference} $\ell$ is a circle of circumference
$1/\ell$.}

Let $Y\subset X$ be the open orbit, equipped with
the K\"{a}hler structure inherited from $X$. Let
$\tY\to Y$ be the universal cover, equipped with the pull
back K\"{a}hler structure. The complex and symplectic
structures on $\tY$ and $Y$ are summarized in
the following diagram:

$$
\begin{CD}
TN_\bR = &N_\bC&  \stackrel{\textup{complex}}{\cong} & \tY&  \stackrel{\textup{anti-symplectic}}{\cong}
& T^*\tri^\circ&  \subset T^*M_\bR\\
& @VV{\exp_T}V   @VVV  @VVV\\
& T&  \stackrel{\textup{complex}}{\cong} & Y&  \stackrel{\textup{anti-symplectic}}{\cong}
& T_\bR\times \tri^\circ \\
& @VV{\log_T}V  @VV{id_Y}V  @VV{\mu_{T_\bR}}V \\
&N_\bR @<{\pi}<< Y @>{p}>> \tri^\circ &\subset M_\bR
\end{CD}
$$
In the above diagram,
$$
\exp_T: N_\bC \cong \bC^n \to T\cong (\bC^*)^n,\quad
(w_1,\ldots, w_n)\mapsto(e^{w_1},\ldots, e^{w_n}),
$$
is the exponential map from the Lie algebra of $T$ to the Lie group $T$,
and
$$
\log_T: T\cong (\bC^*)^n \to N_\bR \cong\bR^n ,\quad
(t_1,\ldots,t_n)\mapsto (\log|t_1|,\ldots, \log|t_n|),
$$
is the logarithm map.
$T^* \tri^\circ\cong N_\bR \times \tri^\circ$ is equipped with
the canonical symplectic form
$$
d(\sum_{i=1}^n   \theta_i dx_i) = \sum_{i=1}^n d\theta_i\wedge dx_i.
$$
This descends to a symplectic form on $T_\bR\times \tri^\circ = (N_\bR/N)\times \tri^\circ$.
The map $T^*\tri^\circ \cong N_\bR \times \tri^\circ \to T_\bR\times \tri^\circ$
is given by $(\theta,x)\mapsto (\exp_{T_\bR}(\theta),x)$. Here
$$
\exp_{T_\bR}: N_\bR\cong \bR^n \to T_\bR\cong U(1)^n,
\quad (\theta_1,\ldots, \theta_n)\mapsto
(e^{\sqrt{-1}\theta_1},\ldots, e^{\sqrt{-1}\theta_n})
$$
is the exponential map from the Lie algebra $\ft_\bR$ to
the Lie group $T_\bR$; the kernel $N\subset N_\bR$
of $\exp_{T_\bR}$ is given by $\theta_i \in 2\pi\bZ$.
The map $T_\bR\times \tri^\circ \to \tri^\circ$ is
the projection to the second factor; it is also the moment map of the
$T_\bR$-action on $T_\bR \times \tri^\circ$ by
multiplication on the first factor. Note that
\cite[Theorem 6.4]{Gu}
$$
\int_X \frac{\omega^n}{n!}=
\int_Y\frac{(\sum_{i=1}^n d\theta_i \wedge dx_i)^n}{n!} =
(2\pi)^n \int_{\tri^\circ} dx_1 \cdots  dx_n
=(2\pi)^n \mathrm{volume}(\tri).
$$

Following \cite{LYZ},
one may apply $T$-duality to the $T_\bR$-fibrations $\pi:Y\to N_\bR$
and $p: Y \to \tri^\circ$ to obtain
the $T$-dual $Y^\vee$ together with a K\"{a}hler
structure.  This construction
is described explicitly by Auroux \cite[Section 4]{Au1} and
Chan-Leung (\cite[Section 3]{CL1} and \cite[Section 2.1]{CL2}).
In particular, they found that as a complex manifold, $Y^\vee$ is an open
subset of the Hori-Vafa mirror $T^\vee \cong (\bC^*)^n$.

Let $\tY^\vee\to Y^\vee$ be the universal cover, equipped with the pull
back K\"{a}hler structure. The complex and symplectic
structures on $\tY^\vee$ and $Y^\vee$ are summarized in
the following diagram:
$$
\begin{CD}
T^*N_\bR&  \stackrel{\textup{anti-symplectic}}{\cong} & \tY^\vee&  \stackrel{\textup{complex}}{\cong}
& M_\bR \times \tri^\circ&  \subset &M_\bC &= TM_\bR\\
 @VVV   @VVV  @VVV @VV{\exp_{T^\vee}}V\\
 T^\vee_\bR\times N_\bR &  \stackrel{\textup{anti-symplectic}}{\cong} & Y^\vee&  \stackrel{\textup{complex}}{\cong}
& \log_{T^\vee}^{-1}(\tri^\circ) &\subset&  T^\vee @>\textup{holomorphic}>W>\bC\\
 @VV{\mu_{T_\bR^\vee}}V  @VV{id_{Y^\vee} }V  @VVV @VV{\log_{T^\vee}}V \\
N_\bR @<{\pi^\vee}<< Y^\vee @>{p^\vee}>> \tri^\circ &\subset& M_\bR
\end{CD}
$$
In the above diagram,
$\exp_{T^\vee}: M_\bC \cong \bC^n \to T^\vee\cong (\bC^*)^n$
is the exponential map from the Lie algebra of $T^\vee$ to the Lie group
$T^\vee$, and
$\log_{T^\vee}: T\cong (\bC^*)^n \to N_\bR \cong\bR^n$
is the logarithm map. $T^*N_\bR \cong M_\bR \times N_\bR$ is equipped with
the canonical symplectic form
$$
d(\sum_{i=1}^n   \gamma_i dy_i) = \sum_{i=1}^n d\gamma_i\wedge dy_i.
$$
This descends to a symplectic form on $T_\bR^\vee\times N_\bR  = (M_\bR/M)\times N_\bR$.
The map $M_\bR \times N_\bR \to T_\bR^\vee \times N_\bR$
is given by $(\gamma,y)\mapsto (\exp_{T_\bR^\vee}(\gamma),y)$.
Here $\exp_{T_\bR^\vee}: M_\bR \to T^\vee_\bR$
is the exponential map from the Lie algebra $\ft_\bR^*$ to the Lie group $T^\vee_\bR$;
the kernel $M\subset M_\bR$ is given by $\gamma_i\in 2\pi\bZ$.
The map $T_\bR^\vee\times N_\bR \to N_\bR$ is
the projection to the second factor; it is also the moment map of the
$T_\bR^\vee$-action on $T_\bR^\vee \times N_\bR$ by
multiplication on the first factor. Finally, we have
$$
\begin{CD}
T^*M_\bR &\stackrel{\textup{symplectic}}{\cong} & \tY^\vee\\
@VVV @VVV\\
T^* T_\bR^\vee &\stackrel{\textup{symplectic}}{\cong} & Y^\vee
\end{CD}
$$
where $T^* M_\bR$ and $T^* T_\bR^\vee$ are equipped with the canonical
symplectic form of the cotangent bundle:
$$
d(\sum_{i=1}^n y_i d\gamma_i) = \sum_{i=1}^n dy_i\wedge d\gamma_i.
$$

\begin{example}
$X= \bP^1$.

The moment map $\mu_{T_\bR}: \bP^1\to M_\bR \cong \bR$ is given by
$$
\mu_{T_\bR}(X_0,X_1)=\frac{a |X_1|^2}{|X_0|^2 + |X_1|^2},\quad a>0.
$$
The image is $[0,a]\subset \bR$.

The complex coordinates on $\tY\cong \bC$ and $Y\cong \bC^*$ are
$y+\sqrt{-1}\theta$ and $e^{y+\sqrt{-1}\theta}=X_1/X_0$, respectively.
The restrictions of the symplectic form $\omega$ and the Riemannian metric $g$
on $X=\bP^1$ to $Y=\bC^*$ are given by
\begin{eqnarray*}
\omega &=& dx \wedge d\theta
= \frac{2a e^{2y} dy\wedge d\theta}{(1+ e^{2y})^2},\quad y\in N_\bR \\
g &=&  \frac{a}{2x (a-x)} dx^2 + \frac{2x(a-x)}{a} d\theta^2
=\frac{2a e^{2y}}{(1+e^{2y})^2}(dy^2 +d\theta^2),\quad y\in \bR
\end{eqnarray*}
where
$$
x\in (0,a)\subset M_\bR \cong \bR,\quad
y \in N_\bR\cong \bR,\quad \theta \in N_\bR/N \cong \bR/2\pi\bZ.
$$

The symplectic form $\omega^\vee$ and the Riemannian metric $g^\vee$
on $Y^\vee$ are given by
\begin{eqnarray*}
\omega^\vee &=& \frac{a dx \wedge d\gamma}{2x(a-x)} = dy\wedge d\gamma,\\
g^\vee &=&  \frac{a}{2x (a-x)} (dx^2 + d\gamma^2)
=\frac{2a e^{2y}}{(1+e^{2y})^2}dy^2 +\frac{(1+e^{2y})^2}{2a e^{2y} }d\gamma^2.
\end{eqnarray*}
where
$$
x\in (0,a)\subset M_\bR \cong \bR,\quad
y \in N_\bR\cong \bR,\quad \gamma \in M_\bR/M\cong \bR/2\pi\bZ.
$$

Recall that under the T-duality, a circle of radius $R$
is dual to a circle of radius $1/R$. In this example,
the fiber $p^{-1}(x)$ of $p:Y\to \tri^\circ=(0,a)$ is a circle of radius
$\displaystyle{\sqrt{\frac{2x(a-x)}{a}}}$, while
the fiber $(p^\vee)^{-1}(x)$ of $p^\vee:Y^\vee\to \tri^\circ=(0,a)$ is a circle of radius
$\displaystyle{\sqrt{\frac{a}{2x(a-x)}}}$.

\end{example}

\subsection{Statement of symplectic results}
\label{subsec:statement-fltz}
We now describe the symplectic results from \cite{FLTZ}.  The description
of the T-duality transformation of
objects is postponed until Section \ref{subsec:tobj}

A smooth projective toric variety is defined by
a fan $\Sigma \subset N_\bR$. Fans of some toric surfaces are shown
in Figure 3, where $\bF_m$ is the Hirzebruch surface
$\bP\left(\cO_{\bP^1}\oplus \cO_{\bP^1}(m)\right)$ $(m\geq 0)$,
and $B_k$ is the toric blowup of $\bP^2$ at $k$ points.
\begin{figure}[h]
\psfrag{P2}{$\bP^2$}
\psfrag{P1P1}{$\bF_0=\bP^1\times \bP^1$}
\psfrag{dP1}{$B_1=\bF_1$}
\psfrag{dP2}{$B_2$}
\psfrag{dP3}{$B_3$}
\psfrag{F2}{$\bF_2$}
\psfrag{F3}{$\bF_3$}
\begin{center}
\includegraphics[scale=0.6]{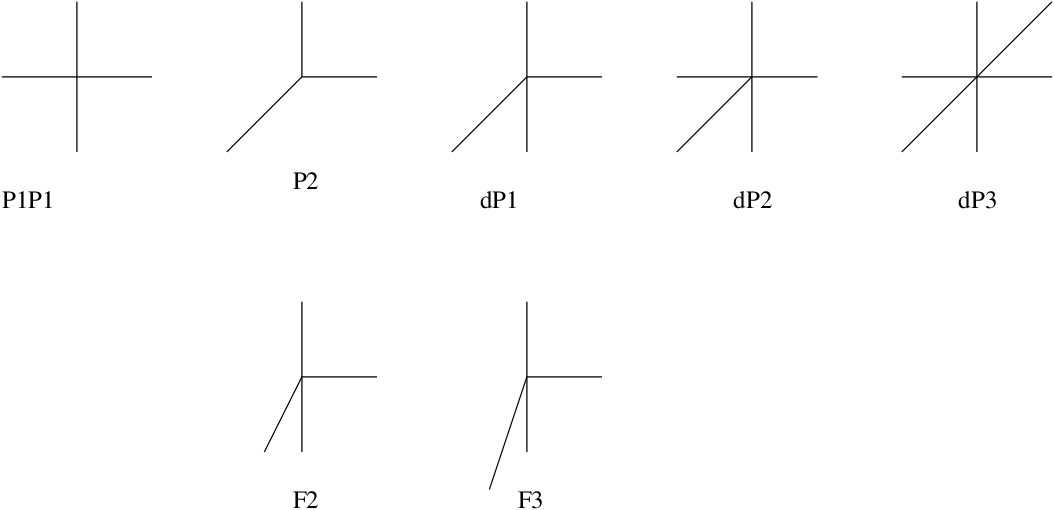}
\end{center}
\caption{Fans of some toric surfaces}
\end{figure}

The first row in Figure 3 consists of all projective smooth toric
Fano surfaces.
Note that the anti-canonical divisor $-K_{\bF_m}$ is
ample for $m=0,1$, numerically effective for $m=2$, and
not numericaly effective for $m>2$. In general, the fan of a smooth
projective toric surface is determined by its 1-dimensional cones.

Let $\langle \   , \ \rangle: M_\bR\times N_\bR \to \bR$ denote
the natural pairing. Given a $d$-dimensional cone
$\sigma \subset N_\bR$, define an $(n-d)$-dimensional
subspace
$$
\sigma^\perp=\{ m\in M_\bR \mid \langle m,v\rangle = 0\  \forall v\in \sigma\}
\subset M_\bR.
$$

Let $X_\Si$ denote the toric variety defined by
a fan $\Si$. Given a fan, define a conical Lagrangian
submanifold of $M_\bR \times N_\bR \cong T^*M_\bR$:
\begin{equation}\label{eqn:LS-M}
\LS = \bigcup_{\sigma\in \Si} (\sigma^\perp + M)\times (-\sigma)
\end{equation}
where $\sigma^\perp + M=\{ x + m\mid x\in \sigma^\perp,\ m\in M\}$ is
the union of all translations of $\sigma^\perp$ by a point in
the lattice $M$. Note that $\LS$ descends to a conical Lagrangian
\begin{equation}\label{eqn:LS-T}
\LS/ M \subset (M_\bR/M)\times N_\bR \cong T^* T_\bR^\vee.
\end{equation}

For example, when $\dim_\bC X_\Si=2$ is a projective surface,
$\Lambda_\Si$ is determined by the 1-dimensional
cones $\rho_i$, $i=1,\ldots,r$. Each 1-dimensional
cone $\rho_i$ determines a circle
$$
\rho_i^\perp/ M \subset M_\bR/ M \cong (S^1)^2.
$$
The fiber of $\LS/ M \to M_\bR/  M$ over
a (generic) point in $\rho_i^\perp/ M$ is
$-\rho_i \subset N_\bR$.
We use an inner product to identify $N_\bR$ with
$M_\bR$ and view $-\rho_i$ as a ray normal to
$\rho_i^\perp/ M$. Therefore
$\LS$ can be drawn on a square:

\begin{figure}[h]
\psfrag{P2}{$\bP^2$}
\psfrag{F0}{$\bF_0=\bP^1\times \bP^1$}
\psfrag{F1}{$B_1=\bF_1$}
\psfrag{dP2}{$B_2$}
\psfrag{dP3}{$B_3$}
\psfrag{F2}{$\bF_2$}
\psfrag{F3}{$\bF_3$}
\begin{center}
\includegraphics[scale=0.6]{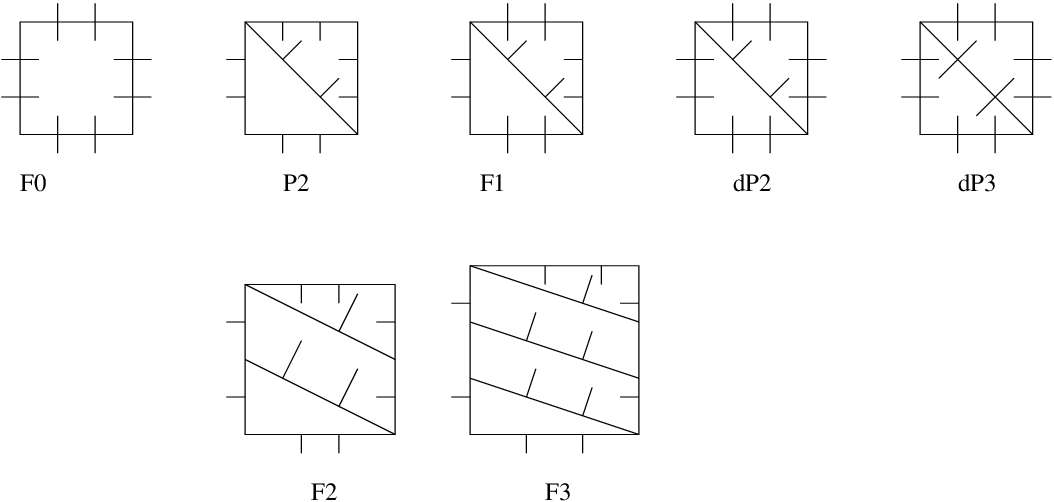}
\caption{Conical Lagrangians of some toric surfaces}
\end{center}
\end{figure}

Let $\iota: T^*M_\bR \hookrightarrow D^*M_\bR$
be the diffeomorphism from the cotengent bundle
to its disk bundle. Given a Lagrangian submanifold or
Lagrangian subvariety $L\subset T^*M_\bR$,
defined the conic limit of $L$ to be
$$
L^\infty = \overline{\iota(L)}\cap S^*M_\bR \subset S^*M_\bR,
$$
where $S^*M_\bR$ is the sphere bundle of the cotangent
bundle $T^*M_\bR$. The conic limit of a Lagrangian
submanifold of $T^*T^\vee_\bR$ is defined similarly.

In \cite{FLTZ}, the authors
derived the following equivariant version of homological mirror symmetry (``(b)'')
for any projective toric manifold:
\begin{eqnarray*}
DCoh_T(X_\Si)&\cong& D Fuk(T^*M_\bR;\Lambda_\Si).
\end{eqnarray*}
Here $Coh_T(X_\Si)$ is the category of $T$-equivariant
coherent sheaves. $D Coh_T(X_\Si)$ is generated
by $T$-equivariant ample line bundles.
An object in $Fuk(T^*M_\bR;\LS)$ is an
exact Lagrangian submanifold $L$ of $T^*M_\bR$
with compact horizontal support such that
$L^\infty\subset \Lambda_\Si^\infty$.
The hom spaces are defined as in Section \ref{sec:fukaya}.

\subsection{T-dual of an equivariant line bundle}
\label{subsec:tobj}
The T-duality functor expressing the equivalence above is constructed
on a generating set of holomorphic objects, equivariant ample line
bundles.  We describe this here.

Let $D_i$, $i=1,\ldots,r$ be
the codimension one orbit closures of $X$.
Then any $T$-invariant divisor is of the form
$$
D_\vc = \sum_{i=1}^r c_i D_i,\quad  \vc=(c_1,\ldots,c_r)\in \bZ^r,
$$
and  any $T$-equivariant line bundle is of the form $\cO_X(D_\vc)$.
There is a $T$-invariant meromorphic section $s$ of $\cO_X(D_\vc)$,
unique up to multiplication by $\bC^*$, such that
$\mathrm{div}(s)=D_\vc$. $s$ is holomorphic and nonvanishing on $Y$,
so it is a holomorphic frame of $\cO_X(D_\vc)\bigr|_Y$. Let $\nabla$
be the $U(1)$-connection determined by a $T_\bR$-invariant hermitian
metric $h$ on $\cO_X(D_\vc)$, and let $\alpha$ be the connection
1-form with respect to the unitary frame $s/\| s \|_h$. The
connection $1$-form $\alpha$ is an \emph{equivariant connection}, as
discussed in \cite{AB}. Then $\alpha$ is purely imaginary, and the
restriction of $\sqrt{-1}\alpha$ to a fiber $Y_x : = p^{-1}(x)$ of
the $T_\bR$-fibration  $p:Y\to \tri^\circ$ is a  real harmonic
1-form on $Y_x$, which can be identified with an element in
$H^1(Y_x;\bR)$, the universal cover of $Y^\vee_x \cong H^1(Y_x;\bR)/
H^1(Y_x;\bZ)$. Letting $x$ vary yields a section of $\tilde{Y}^\vee
\cong M_\bR\times N_\bR \to N_\bR$, which is the same as a map
$\Psi_h: N_\bR \to M_\bR$. The map $\Psi_h$ has the following
interpretation. Let $F_h = d\alpha$ be the curvature 2-form. Then
$\omega_h=\sqrt{-1}F_h$ is a closed $T$-invariant 2-form (a {\em
presymplectic form} in the sense of Karshon-Tolman \cite{KT}). Note
that
$$
[\omega_h] = 2\pi c_1(\cO_X(D_\vc)) \in H^2(X;\bR).
$$
$\omega_h$ defines a moment map $\Phi_h: X\to M_\bR$ up to a
constant vector in $M_\bR$, given explicitly in \cite{AB}; the
constant is determined by the equivariant structure on
$\cO_X(D_\vc)$. Let $j_0: N_\bR \to X$ be the composition of $\exp:
N_\bR \to Y$ and the open embedding $Y\to X$. Then $\Psi_h=
\Phi_h\circ j_0$. The T-dual Lagrangian of the the $T$-equivariant
line bundle $\cO_X(D_\vc)$ equipped with a $T_\bR$-invariant
hermitian metric $h$ is given by
$$
\bL_{\vc,h}=\{ (\Phi_h\circ j_0(y),y)\mid y\in N_\bR\}  \subset M_\bR\times N_\bR = T^*M_\bR.
$$

\begin{example}
For $\bP^1$, the $1$-cones in the fan $\Si$ are generated by $v_1=1,\ v_2=-1$. Let $D_1$ and $D_2$ be two $T$-invariant divisors corresponding to these $1$-cones, and
$z$ be the complex coordinate such that $z|_{D_1}=0$ and $z|_{D_2}=\infty$. The exponential map $j_0: N_\bR\rightarrow \bP^1$ is simply given by $y_1\mapsto e^{y_1}$. The $T$-equivariant line bundle $\cO(c_1 D_1+ c_2D_2)$ has a canonical meromorphic section $s$. Endow this line bundle  with a hermitian metric $h$
$$\|s\|_h^2=\frac{|z|^{2c_1}}{(1+|z|^2)^{c_1+c_2}}.$$
The moment map $\Phi_h$ associated to the presymplectic form $\omega_h=\sqrt{-1}F_h$ is given by
$$z\mapsto \frac{(c_1+c_2)|z|^2}{1+|z|^2}-c_1.$$
Thus the map $\Psi_h: N_\bR\to M_\bR$ is given by the formula
$$
\frac{\gamma_1}{2\pi} =
\frac{(c_1+c_2)e^{2y_1}}{1+e^{2y_1}}-c_1.
$$

This equation characterizes the T-dual Lagrangian $\bL_{\vc,h}$ where $\vc=c_1v_1+c_2v_2$.

\end{example}

$D Coh_T(X)$ is generated by equivariant ample line bundles,
so it suffices to define the functor
$D Coh_T(X)\to D Fuk(T^*M_\bR;\LS)$ on equivariant
ample line bundles.
When $\cO_X(D_\vc)$ is an ample, there exists a $T_\bR$-invariant
hermitian metric $h$ on $\cO_X(D_\vc)$ such that
$\omega_h$ is a symplectic form. Then the image of
the moment map $\Phi_h$ is a convex polytope
$$
\tri_\vc =\{ m\in M_\bR\mid \langle m, v_i\rangle\geq -c_i\mid i=1,\ldots,r\},
$$
where $v_i\in N$ is primitive, and $\bR_{\geq 0} v_i$ is the
1-dimensional cone associated to $D_i$.
The map $\Psi_h =\Phi_h\circ j_0: N_\bR\to \tri_\vc^\circ$ is
a diffeomorphism from $N_\bR$ to the interior of $\tri_\vc$.
Let $F_\sigma$ be the face of $\tri_\vc$ associated to
the cone $\sigma\in \Sigma$, and define a conical Lagrangian
$$
\Lambda_\vc = \bigcup_{\si\in \Sigma} F_\sigma\times(-\sigma) \subset \LS.
$$
Then $\bL_{\vc,h}^\infty = \Lambda_\vc^\infty\subset \LS^\infty$. It is
indeed an object in $Fuk(T^*M_\bR;\LS)$. Different metrics
define equivalent objects.

\section{Microlocalization}
\label{sec:micro}

We are proving an equivalence between
coherent sheaves, Lagrangians, and constructible sheaves.
The first equivalence is described in the preceding section.
The last equivalence -- between the Fukaya category $Fuk(T^*B)$
of a cotangent bundle
and constructible sheaves $Sh_c(B)$ on the base -- is provided by the papers \cite{NZ,N},
whose results and ideas we describe in this section.
The strategy in relating the Fukaya category of a cotangent bundle to
constructible sheaves will be to exploit how both are computed
via Morse theory.  The relation to Morse theory is found
following the lead of
Kontsevich-Soibelman \cite{KoSo} and
Fukaya-Oh \cite{FO}, who studied a similar
geometry.

\subsection{Algebraic preliminaries}\label{sec:algebra}
Here we review some of the algebraic structures which appear
in the argument.
We will work with dg and $A_\infty$ categories.
Recall that a dg algebra is a complex with differential $d$
and a degree-zero product $\circ$ obeying the Leibnitz rule.
The path from dg algebras to dg categories
is clear:  the
the hom spaces in dg categories have
the structure of a differential complex, with
compositions morphisms of complexes.
$A_\infty$ algebras have degree-one differentials $m_1 := d$,
Leibnitz degree-zero products $m_2 = \circ,$ and higher compositions $m_k$
of degree $2-k$ (so a dg algebra is an $A_\infty$ algebra with $m_{k\geq 3}=0$).
Recall here the central example
of the $A_\infty$ algebra of chains on a loop
space:\footnote{We refer the reader to
\cite{K} for rigorous definitions and an excellent discussion.}
the concatenation product is not associative,
but is associative at the level of homology.  Further, the different
homotopies from $((a*b)*c)*d$ to $a*(b*(c*d))$ are themselves
homotopic, and so on, so there are higher relations as well.
({\em N.B.}:  The suspension from
a space to its loops helps to explain the unusual grading
conventions one often encounters.)
There is a clear path from
$A_\infty$ algebras to $A_\infty$ categories.

Given an $A_\infty$ category $\cC$ we can consider
the dg category $\cC$-mod of (contravariant)
$A_\infty$ functors from $\cC$ to the
category of chain complexes.  The Yoneda embedding
is a functor from $\cC$ to $\cC$-mod sending an object
$a$ to $hom_\cC(-,a).$

If $\cC$ is a dg or $A_\infty$ category (the former
a special case of the latter), then taking cohomology
of the hom spaces yields an ordinary (associative)
$\bZ$-graded category $H(\cC)$.   $H(\cC)$ has
a triangulated structure:  a triangle is distinguished
if its image under Yoneda is (isomorphic to) an exact triangle.
$A_\infty$ categories with shift functors satisfying
some conditions are said to be triangulated.
The triangulated envelope of $\cC$, denoted $Tr\cC$,
is, informally, the smallest $A_\infty$ category generated by
$\cC$.  It is unique up to quasi-equivalence, and
can be constructed explicitly from twisted
complexes or as the envelope within $\cC$-mod
of the image under the Yoneda embedding.

\subsection{The Cast of Categories}
\label{cast}

Some intermediate categories
between $Fuk(T^*B)$ and $Sh_c(B)$
figure prominently in the proof of equivalences.  After
recalling the definition of $Sh_c(B)$ below, we will
describe these
intermediate categories and the quasi-equivalences between
them.

\subsubsection*{$Sh_c(B)$}
We recall what we mean by
the category of constructible sheaves, $Sh_c(B).$
Recall that a constructible sheaf on
a topological space $B$ is a sheaf (of $\bC_B$ modules)
that is locally constant on the strata of a Whitney
stratification.\footnote{We recall that a stratification
is a decomposition of a manifold into strata of pairwise disjoint
manifolds such that one intersects the closure of another
if and only if it is contained therein.  The Whitney condition
avoids some kinds of pathological behavior on the
local behavior of different strata.}
We use the term ``constructible sheaf'' (sometimes
just ``sheaf'') to describe
a complex of sheaves with bounded, constructible
cohomology.  Acyclic complexes form a null system within
this dg category,
and we write $Sh_c(B)$ for the associated quotient.
Despite the quotient, $Sh_c(B)$ retains the structure of a $dg$ category.
Taking
cohomology in degree zero produces $DSh_c(B),$ the
bounded derived category of constructible sheaves on $B.$
As an example of an object, to an open submanifold $U\stackrel{i}{\hookrightarrow}B$
(compatible with a Whitney stratification), we can associate
the sheaf $i_*\bC_U.$  We call this a ``standard open'' sheaf.
In this section, we will assume $B$ is compact.
When $B$ is noncompact we write $Sh_{cc}(B)$ for the subcategory
with compactly supported cohomology\footnote{So the first subscript "c" is
for "constructible," the second for "compact."}
and $DSh_{cc}(B)$ for its derived category.

\subsubsection{The Morse Category}
As discovered by Fukaya \cite{F} (and Cohen-Jones-Segal \cite{CJS}), Morse theory
contains an $A_\infty$ category structure.  Define $Mor(B)$ to be the
category with objects consisting of open sets with defining functions for their
boundaries,
i.e. pairs $(U,m)$ with $U\subset B$
an open submanifold and $m:B\rightarrow \bR$ a function
obeying $m\vert_U >0,$
$m\vert_{\partial U} = 0.$  We put $f = \log m\vert_U : U\rightarrow \bR.$  Hom spaces
are made from the critical points
graded by Morse index, that is
$hom_{Mor(B)}((U,m),(U',m')) = \bigoplus_{p\in Crit(f'-f)}\bC_p[-deg(p)]$.  The differential $d$ is
defined by counting (with sign) points in the space of gradient flow trajectories
$\{\dot{\gamma} = \nabla (f'-f)\}$ modulo time translation.
Higher compositions, such as the product,
are defined by counting gradient flow {\em trees}.
For example, if we want to compute a product $p\circ q$
and find the coefficient of a critical point $r,$ we count
(with sign) trajectories of gradient flow lines for $f'-f$
emanating from $p$ that meet flows for $f''-f'$ from $q$
(at some point) and continue flowing by $f''-f$ to $r$.
The ``tree'' here is shaped like the letter {\sf Y}.\footnote{In case it
was not clear from the context, $p\in hom((U,m),(U',m')),$
$q\in hom((U',m'),(U'',m''))$ and $r\in hom((U,m),(U'',m'')).$}

The conditions on $m$ ensure that gradient flows of $f = \log m$ do
not leave the open sets, and this local version of Morse theory
makes sense.

\subsubsection{The De Rham Category}

In fact, this $A_\infty$ structure is a consequence of
the $dg$ structure of a de Rham category we now describe.
The objects of $DR(B)$ are also pairs $(U,m)$
and the
hom spaces are relative de Rham
complexes
$hom_{DR(B)}(U,U') = (\,\Omega^*(\overline{U}\cap U'),\partial{U}\cap U'\, ,d\,)$
(if one likes, one can omit the inconsequential data of $m$ in the
objects).  For example, when $U=U'=B$ we get
$\Omega^*(B)$ with de Rham differential.

The point is that a $dg$ algebra $A$ induces an $A_\infty$
structure on a subvector space $B$ which appears as a
strong deformation retract, i.e. the image of a projector.
That is, we have $p:A\rightarrow A,$ $B = {\rm Image}(p),$
$i:B\hookrightarrow A$ with $p \circ i = id_B$ and a homotopy $h:A\rightarrow B$ of degree $-1$
such that $dh + hd = i \circ p - id_A.$
A simple example is when $A$ is the de Rham complex of a Riemannian
manifold and $B$ is the vector space of harmonic differential forms.
Let $\triangle$ be the Laplacian and $G$ the Green operator.
Then $i$ is inclusion, $p$ is projection ${\bf 1} - G\triangle$ and
$h = Gd^*.$

The induced compositions $m_k(b_1,...,b_k)$ are defined as follows.
First define $a_j = i(b_j)$ from the inclusion, then form a
sum over all trivalent
trees connecting vertices $a_j$ to a terminal root, where each tree
counts by forming the dg product $\circ$ from $A$ at each internal vertex
and the homotopy $h$ is applied at each internal edge.  Finally, the
projection $p$ is applied at the terminal root so that the result
lies in $B$. For example, the contribution to
$m_4(b_1,b_2,b_3,b_4)$ from the tree ${}^{\vee\vee}\!\!\!\!\!{}_{\textsf Y}$
is $p(h(i(b_1)\circ i(b_2))\circ h(i(b_3)\circ i(b_4))).$
Then $i$ and $p$
are quasi-equivalences of $A_\infty$ algebras.
The same construction holds
for categories.

Harvey and Lawson found the homotopy relating $DR(B)$
to $Mor(B)$ using the language of currents.  As any submanifold
of $B\times B$ or ``kernel'' defines, through pull-back
and push-forward (integration), an operator
on currents (which we will describe as differential
forms for simplicity).  Given a function $f,$ the kernel corresponding
to gradient flow $\varphi_t$
for time $t$ is the graph $\Gamma_t :=
\{x,\varphi_t(x)\}$.
Letting $t$ range from $0$ to a fixed point at $t=\infty,$ we get a homotopy
between the identity $\varphi_0$ and the space of currents supported
at critical points (thought of as the Morse complex vector space),
given by the union:  $h =
\cup_{0\leq t\leq \infty}\Gamma_t.$

\subsubsection{The Category $Open(B)$}
\label{open}

The category $DR(B)$ looks a lot like the dg subcategory $Open(B)\subset Sh_c(B)$
whose objects are again pairs $(U,m)$, thought of as representing
the ``standard open'' sheaves $i_*\bC_U$ (again, $m$ is inconsequential).
In fact, one can easily show
\begin{equation}
\label{comps}
hom_{DR(B)}(U,U') = (\,\Omega^*(\overline{U}\cap U'),\partial{U}\cap U'\, ,d\,)\cong hom_{Sh_c(B)}(i_*\bC_U,i_*\bC_{U'}).
\end{equation}
Further, the compositions coincide.
As a result, these categories have the same triangulated
envelope.  But the triangulated envelope
$TrOpen(B)$ is all of $Sh_c(B)$!  Informally,
``standard opens generate.''\footnote{This follows from two
facts.  First, any complex of sheaves can be constructed froms
cones and shifts of
its cohomology sheaves, so it remains to prove that constructible
sheaves are generated by standard opens.  This follows by triangulating
the space and writing a non-open simplex $T$ as the ``difference'' between
the open ``star'' $Star(T)$ of all simplices containing $T$ in their closure
and the complement $Star(T)\setminus T,$ with $T$ removed.}

Taken together, these results show the following
$A_\infty$ quasi-equivalences:
$$TrMor(B) \cong TrDR(B) \cong Sh_c(B).$$
What remains is to show that $Mor(B)$ is equivalent
to the Fukaya category $Fuk(T^*B)$.  To do so, first, we find a
quasi-embedding,
then we construct an inverse.

\subsection{Fukaya-Oh Theorem}
\label{fothm}

Central to the construction of an equivalence
is the theorem of Fukaya and Oh \cite{FO} relating
the subcategory of $Mor(B)$ defined by global functions, i.e. objects $(U=B,m=\exp(f))$,
to the subcategory of $Fuk(T^*B)$ defined by global graphs $\Gamma_{df}.$
To explain the idea, let us
recall that Fukaya theory is like Morse theory on the space of paths $c$
between pairs of Lagrangian submanifolds (for the Morse action functional
$a(c) = \int_c \theta,$ where $\theta$ is the canonical one-form on
$T^*B).$\footnote{Let $\pi:T^*B\rightarrow B$ be the projection.  Then $\theta_{(x,\xi)}(v)
= \xi(\pi_* v).$  On a general symplectic manifold we put $a(c) = \int_D \omega,$ where
$D$ is a homotopy between $c$ and some fiducial path $c_0$.}
Floer \cite{Fl} related this Morse theory for pairs of Lagrangian graphs in a cotangent bundle
to ordinary Morse theory on the base.
Specifically, he showed that pseudoholmorphic strips between
such Lagrangian graphs are in correspondence
with gradient flow lines for the difference of the Morse functions.
Fukaya and Oh extended this idea to more general
pseudoholomorphic disks.\footnote{In a recent preprint of Iacovino \cite{I},
a cleverly chosen inhomogeneous term is added to the pseudoholomorphc
curve equation to
show that the image of such a pseudholomorphic disk {\em is}
a gradient flow tree.}

Consider a pair of functions $(f,f')$ and a pair of
exact Lagrangian graphs $\Gamma_{\epsilon
df},\Gamma_{\epsilon df'}.$
Fukaya and Oh proved not only that for small enough $\epsilon$
the Morse moduli space of gradient flow lines for $f'-f$ is diffeomorphic to the
moduli space of holomorphic strips bounding the
two Lagrangians.  They also showed that the same is true for
collections $\vec{f} = (f_0,...,f_k)$ and $\Gamma_{\epsilon
d\vec{f}}=
(\Gamma_{\epsilon df_0},...,\Gamma_{\epsilon df_k})$
when we consider the Morse moduli space of gradient flow
trees and the Fukaya moduli space of disks.  The case of
three objects with a gradient flow tree a letter {\sf Y} is illustrative:
the corresponding disk looks like a thickening.  Fukaya and Oh
construct from the {\sf Y} tree an approximate holomorphic disk,
then prove that an actual holomorphic disk exists nearby.
The ``nearby'' aspect is important, since it is clear from the discussion
in Section \ref{cast} that we will need a local version of the Fukaya-Oh result.

\subsection{Building the Equivalence}

This discussion, plus our notations, suggests that to construct an
embedding from $Mor(B)$ to $Fuk(T^*B)$ we should map $(U,m)$
to $\Gamma_{df}$, a Lagrangian graph over $U$.
In fact, as discussed in Section \ref{open},
the category $Sh_c(B)$ is generated by ``standard open''
sheaves $i_*\bC_U$, so we can
define a functor by defining it on the full subcategory $Open(B)$ consisting
of these objects and
extending it as a triangulated functor to
the triangulated envelope, which is $Sh_c(B)$.  This
will be our strategy for proving $Sh_c(B) \cong Fuk(T^*B).$

Two simple examples should help develop some intuition
for why the embedding is constructed this way.

\begin{example}

Consider the constant sheaf
$\bC_B.$  Then $hom_{Sh_c(B)}(\bC_B,\bC_B)$
can be taken to be any model for sheaf cohomology\footnote{Recall that in $Sh_c(B)$
we have quotiented by homotopies.} $H^*(B,\bC_B)\cong H^*(B)$ such as
$\check{\rm C}$ech or Morse complexes
(or de Rham, if we take coefficients in $\bR$ or $\bC$).
The Morse complex is obtained by first choosing a Morse function
(which we think of as a section of $\bC_B\otimes_{\bC_B}C^\infty(B)$)
and then analyzing its gradient flow trajectories.

We now recall that
the Piukhinin-Salamon-Schwartz (PSS) isomorphism in Floer theory
relates the Floer cohomology of a Lagrangian $L\subset M$ (which we think
of as $H^* hom_{Fuk(M)}(L,L)$) to the
ordinary (e.g., singular) cohomology $H^*(L)$.  We take this as
the first suggestion of a functor relating
$\bC_B$ to the zero section $B\subset T^*B.$
(The zero section is also
the characteristic cycle of $\bC_B$ -- see Section \ref{subsec:cc}.)

In fact, the Morse function $f$ on $L$
can be thought of as perturbation data for the Lagrangian,
since $\omega^{-1}(df,-)$ is a normal vector field to $L$.
As discussed in Section \ref{fothm} above,
Fukaya and Oh showed that for small enough $\epsilon$
the Morse complex of $\epsilon f$ is equal to the Floer complex
of $L$ and its perturbation.  They further proved that
every moduli space used in computing compositions
for the $A_\infty$ Morse category of functions (gradient flow trees)
$f:B\rightarrow \bR$ is oriented-diffeomorphic
to the Fukaya moduli spaces (holomorphic disks)
for the category of graphs $\Gamma_{df}$ in
the cotangent bundle $T^*B.$

\end{example}

\begin{example}
Now consider a hom from $\bC_B$ to a standard $i_*\bC_U.$
Again, the hom in $Sh_c(B)$ is any model for the sheaf
cohomology of $i_*\bC_U,$ i.e. $H^*(U),$ such as the
de Rham complex.  Morse theory on $U$ would work,
with a proper Morse function $f$.  This can be achieved
by taking $f = \log m,$ where $m>0$ on $U$ and $m=0$
on $\partial U$.

From these observations, we can hope that
a ``local'' version of the Fukaya-Oh theorem is true:
the moduli spaces of holomorphic strips computing
$hom_{Fuk(T^*B)}(B,\Gamma_{df})$
are diffeomorphic to the spaces of gradient flow lines
of $f$.  As shown in \cite{NZ}, this turns out to be
the case.  Even more:  the same is also true for the
moduli space of a composition
among $\Gamma_{df_1},...,\Gamma_{df_k}$
and the corresponding moduli space
for Morse functions $f_1,...,f_k.$
\end{example}

With these motivations, we define a
functor $\mu_0: Open(B) \rightarrow Fuk(T^*B)$ by
$\mu_0((U,m)) = \Gamma_{df},$ where $f = \log m$ as
usual.  Since $TrOpen(B)\cong Sh_c(B),$ $\mu_0$ extends
to a {\em microlocalization} functor
$\mu: Sh_c(B)\rightarrow Fuk(T^*B)$
sending $i_*\bC_U$ to $\Gamma_{df}$.
By Equation \eqref{comps} and the local Fukaya-Oh theorem,
$\mu$ is a quasi-embedding (an isomorphism on the cohomology
of the hom complexes).
The Verdier dual of direct image $i_*$ is the proper
direct image $i_!$.  It turns out that Verdier duality in $Fuk(T^*B)$
is multiplication by $(-1)$ in the fibers.
We write $L_{U*} = \mu(i_*\bC_U)$ for the ``standard Lagrangian brane''
or simply ``standard Lagrangian''
and $L_{U!} = \mu(i_!\bC_U)$ for the ``costandard brane''
or ``costandard Lagrangian.''  One can extend the definition of standard
branes to non-open submanifolds, as well (we use this below).

\subsection{Equivalence and the Inverse Functor}

The microlocalization functor $\mu: Sh_c(B)\rightarrow Fuk(T^*B)$ is
a quasi-embedding, and we can define an obvious candidate for
its inverse.  We want to associate a sheaf to a Lagrangian.
Recall that T-duality associates a bundle to a section, and a higher
rank bundle to a multi-section.  Here we have an
analogous story, and the key is to note that the stalk
of the sheaf at a point where it is locally constant is generated by the points
of the Lagrangian over that point.
There may be disks relating the points, however, and
a more invariant definition would involve the whole Fukaya hom {\em complex}
between a fiber and a Lagrangian (thought of as the fiber of a complex
of vector bundles).  We will therefore define for a Lagrangian
a sheaf of complexes. The natural guess (taking contravariance
into account) is
that a Lagrangian brane $P$ is mapped to a sheaf of
complexes $F_P$
defined by
$$F_P^*(U) = hom^*_{Fuk(T^*B)}(L_{U!},P),$$
where $L_{U!}$ is the costandard brane on $U$, as defined above.\footnote{In
fact, this is actually an $A_\infty$ sheaf since the composition of
restrictions will not be associative on the nose -- see \cite{N} for details.}

In \cite{N}, it is proven that the microlocalization functor $\mu$ is a quasi-equivalence
and that the functor $\nu$ mapping $P$ to $F_P$
is a quasi-inverse to $\mu$.  The latter assertion follows readily from the former;
proving the essential surjectivity of $\mu$ is the difficult part.
To do this, Nadler considers the toy problem of showing that
a collection of vectors $v_i$ spans a vector space $V.$  It is enough
to express the identity map in the form ${\bf 1} = \sum v_i \otimes w^i$ for some
(co-)vectors $w^i.$  Then applying $\bf 1$ to $v$ expresses $v$ in terms of the
$v_i.$  Analogously, Beilinson's resolution of the diagonal $\triangle_{\bP^n}\subset
\bP^n\times\bP^n$ (thought of as the kernel for the identity functor) leads
to a resolution of a coherent sheaf in terms of a vector bundles.
A similar trick works here.

First, it suffices to prove surjectivity for
a class of Fukaya
objects that intersect infinity of $T^*B$ inside the boundary-at-infinity of some
fixed (but arbitrary) conical Lagrangian $\Lambda$ at infinity in $T^*B.$  Then
find a triangulation $\{T_\alpha\}$ of $B$ whose associated conical Lagrangian
$\cup_\alpha T^*T_\alpha$ contains $\Lambda.$
Of course the diagonal $\tri_B\subset B\times B$ is not expressible in terms of
products $T_\alpha\times T_\beta$ -- otherwise, we'd nearly be finished.  The trick in
the proof is that the piece of the diagonal $\tri_{T_\alpha} \subset T_\alpha \times T_\alpha$
over $T_\alpha$ is homotopic to $\{t_*\}\times T_\alpha$
for any $t_*\in T_\alpha,$ and Nadler showed that corresponding to this homotopy
is an isomorphism of functors between a piece of the diagonal $Y(L_{\triangle_{T_\alpha}*})$
and an external
product $Y(L_{T_\alpha *}){\tiny \boxtimes} Y(L_{\{t_*\}})$, where $Y$ is the
Yoneda embedding.  Applying the identity operator
in this form shows that the (Yoneda modules corresponding
to) $L_{T_\alpha *}$ generate the category.

\subsection{Singular support and characteristic cycles}
\label{subsec:cc}

One can describe the results relating constructible sheaves to
the Fukaya category
as a categorification of the characteristic cycle
construction
of Kashiwara-Schapira (see \cite{KaSc} for foundational
material).  The path is straightforward,
with some hindsight and a selective look at
some facts involving of characteristic
cycles.
We review this interpretation here.

First recall that given
a constructible complex of sheaves $F$ on $B,$
its {\em singular support}, $SS(F),$ is a conical Lagrangian subvariety in
$T^* B$ which encodes Morse-theoretic obstructions to extending
local sections of $F$.  That is, if a covector is {\em not} in $SS(F)$,
it means that there is no obstruction to propogating local sections of $F$
in directions which are positive on the covector.

There is a finer invariant of $F$ called its {\em characteristic cycle}, or $CC(F)$, which is a linear combination of Lagrangians components of $SS(F)$.  The multiplicity of $CC(F)$ at
a given covector is the Euler characteristic of the local Morse groups of the complex with respect to
the covector (roughly, the restriction map that the sheaf associates to an open neighborhood and to the smaller open set of points evaluating negatively on the covector).

So, for example, the characteristic cycle of a flat vector bundle on $B$ is the zero section in $T^*B$ with multiplicity equal to the rank of the vector bundle.  More generally, the characteristic cycle of a flat vector bundle on a submanifold of $B$ is its conormal bundle with multiplicity equal to the rank of the vector bundle.

\begin{remark}
The singular support is additive on exact triangles.  That is, we have $\SS({\tt Cone}(F \to G)) \subset \SS(F) \cup \SS(G)$ for any morphism of complexes $F \to G$.  It follows that if we fix a conical Lagrangian $\Lambda \subset T^* B$, the full subcategory $\Sh_c(B;\Lambda) \subset \Sh_c(B)$ of sheaves with singular support contained in $\Lambda$ is triangulated.
These subcategories appear repeatedly in applications---most commonly, $\Lambda$ is the conormal variety to a Whitney stratification of $B$, and then $\Sh_c(B;\Lambda)$ is the same as the category of complexes constructible with respect to the
stratification.  The more general notion, when $\Lambda$ is not conormal to a stratification, occurs in our work on toric varieties below.
\end{remark}

The calculation of Schmid-Vilonen \cite{SV} is very suggestive.
These authors computed the characteristic cycle of the
standard open $i_*\bC_U$ to be the limit
$\lim_{\epsilon\rightarrow 0}\Gamma_{\epsilon d\log m}.$
Note that $\epsilon \log m = \log m^\epsilon,$ and as $\epsilon$
becomes small $m^\epsilon$
approaches the indicator function of $U,$
which can be thought of informally as a section of $i_*\bC_U.$
The Morse theory of $\log(m^\epsilon)$ is independent of
$\epsilon>0$.  On the other hand, the limit
 $\epsilon \rightarrow 0$ seems to relate Morse
theory (functions) to constructible
sheaves (constant functions).  In an entirely nonrigorous sense,
the family $m^\epsilon$
is a kind of isotopy between the Morse functions of the Morse
category and the indicator functions of constructible sheaves.

A formula of Kashiwara-Dubson further suggests a relationship
to the Fukaya category.
The K-theory of constructible sheaves is isomorphic to the
abelian group of constructible functions, with the local Euler
characteristic providing the isomorphism.
The characteristic cycle maps this isomorphically to the
group of closed, conical Lagrangian cycles $\mathcal L_{con}(T^*B).$
Dubson-Kashiwara proved that for constructible sheaves $F_1,$ $F_2$
(thought of as elements of the K-theory),
$$\chi(F_1, F_2) = CC({\mathcal D}F_1)\cdot CC(F_2),$$
where $\mathcal D$ is the Verdier duality functor. ($CC(\mathcal DF)=
a(CC(F)),$ where $a$ is (-1) in the fibers of $T^*B$.)
On the left,
the expression
is the descendent in K-theory of the chain complex $hom_{Sh_c(B)}(F_1,F_2)$
computed in the dg category $Sh_c(B).$
On the right, the objects are not taken to lie in any category, but
now we understand the $CC(F)$ as the K-theory descendent of
the brane $\mu(F) \in Ob(TrFuk(T^*B)).$
Therefore,
on the right side we can begin with the complex
$hom_{TrFuk(T^*B)}(\mu(F_1),\mu(F_2))$ and compute
its Euler characteristic to get the required result.
The explicit map from the Fukaya category to its
K-theory (conical Lagrangian cycles)
is simply the limit of dilation toward the zero section
(this leaves infinity intact), which preserves the intersection pairing.

We arrive at the interpretation that the microlocalization
quasi-equivalence is a categorification of the characteristic cycle.

\subsection{Comments on technicalities}

The reader might be dissatisfied with our informal
discussion, and will want to convince her/himself of the details.
We must refer such a reader to the papers \cite{N,NZ};
here, as a preview, we make a few remarks on the technical issues
arising in the proofs, and how they are being dealt with.
(Even still, we omit many.)
Most readers will gladly skip this section.

\begin{enumerate}
\item {\em $Mor(B).$}
To define this category, one must ensure that
the moduli spaces of gradient flows have compactifications, no
matter the behavior of the open sets.  To do so, one uses
stratification theory to conclude that the set $\{m>\epsilon\}$
is a good approximation to $U$ with a smooth boundary, and
similarly for $U'$ and another $\epsilon'$, and that
these two intersect transversely.  One also needs compatibility
with the (possibly perturbed) metric so that the
inward/outward behavior of $\nabla \log f$ and $-\nabla \log f'$ remain
along the ``new'' open sets.
These conditions ensure that sequences of gradient flows
do not ``wander off'' the intersection of the open sets,
so that the standard compactifications by broken trajectories
can be constructed.

\item {\em Brane structures in $Fuk(T^*B).$}
As is beautifully described in \cite{S}, the Fukaya category
requires that $hom$ spaces are graded and that the moduli
spaces defining compositions are oriented (compatibly
with topological field theory gluings).
Also, Lagrangian objects
can carry unitary local systems, which we have not yet discussed.
The existence of a grading follows from a ``grading'' of the
Lagrangian submanifolds, i.e. a lift of the Maslov phase from $U(1)$
to its universal cover $\bR$ (in particular, the Lagrangians must
be Maslov-trivial).  Orientability requires that the Lagrangians
are relatively pin with respect to a background class in $w_2(T^*B)$.
Happily, Lagrangian
graphs over $U$ are homotopic to the zero section along $U$
(which is given the zero grading),
and therefore are canonically graded and relatively pin
with respect to the class $\pi^* w_2(B)$.
By assigning the trivial local system, we find canonical
brane structures on the Lagrangians corresponding to standard
opens.  We call the resulting objects ``standard branes.''

\item {\em Regularity.}
One requires the moduli spaces to be manifolds,
which depends on the vanishing of obstructions to infinitesimal
deformations of each pseudoholomorphic map in the moduli
space.  The situation here is no different than in the compact
Fukaya category (see the excellent \cite{S}) or even symplectic
Gromov-Witten theory.  Regularity is achieved by perturbing the
equations of pseudoholomorphicity.

\item {\em Perturbations and Infinity.}
Nontransversal intersections of Lagrangians
are treated through perturbation data (as in the case of
the compact Fukaya
category).  The novelty in the noncompact case arises from
intersections ``at infinity.''  These are treated by perturbing
with geodesic flow for small time, which is the Hamiltonian perturbation
associated to the distance (from the zero section) function.
For multi-compositions involving several Lagrangians
intersecting at infinity, the composition order determines
the ordering of times of geodesic flow:  perturbations
propagate forward in time:  for example in
$T^*\bR = \bR^2,$ $hom(\{x=0\},\{y=1/x\})= 0$ while
$hom(\{y=1/x\},\{x=0\})\neq 0.$

For these perturbations to separate points at infinity
(and bring Lagrangian intersections into finite space), we
must avoid pathologies such as
intersection points accumulating at infinity (the
helix $\theta = \xi$ in $T^*S^1$ is disallowed).  This leads
us to require that the Lagrangians are ``good'' subsets of
the compactification $\overline{T}^*B,$ such as
subanalytic sets.

\item {\em Compactifying moduli spaces.}
Compositions in $Fuk(T^*B)$ are defined by intersections in moduli spaces
of pseudoholomorphic disks, which therefore must be compactified.  The novelty in the
noncompact case is that sequences of pseudoholomorphic disks
may wander off to infinity.  A tameness condition on the Lagrangians
prevents this from happening \cite{Si}.
In fact, we need a slightly weaker
notion to accommodate standard branes, which are not necessarily tame.
We require simply a tame perturbation -- essentially,
a family of Lagrangians $L_t$ indexed by small $t \geq 0$ which is tame
for all $t> 0$.  Also, an almost complex structure is
chosen on the cotangent bundle which is asymptotically
conical (in particularly, not Sasakian near infinity).
The resulting theory is well-defined.

\item {\em Increasing unions of finite calculations.}  The
arguments here apply to any finite calculations among a finite
collection of branes.  The issue of compatibility for
increasing unions of finite calculations is discussed
in the appendix of \cite{N}.

\item {\em Equivalence and the inverse functor.}
The main analytical issue to deal with here is in expressing the
identity functor in terms of standard branes.  The
diagonal in $B\times B$
is nowhere a Cartesian product.  However, inside the cell $T\times T,$
the diagonal $\tri_T$ is homotopic to $T\times \{t_*\}$ where $t_*$ is any
point in $T.$  Correspondingly, one can construct a family of
Lagrangian branes interpolating between $L_{{\tri_T}*}$
and $L_{T*}\times L_{\{t\}!}$.  The key is that the interpolating
family only moves through infinity without intersecting the
conical set $\Lambda,$ and that this ``noncharacteristic isotopy''
of branes induces isomorphic functors. The analytical argument
proves that the moving moduli spaces of pseudoholomorphic disks which
realize isomorphisms of branes (defined as you would for Hamiltonian
isotopies of compact branes) make sense for small motions.
The trick is a kind of ``shimmy,'' breaking up the motion of a brane
into two parts:  first, the motion near infinity
is done far away from the zero section so as not to affect any
pseudoholomorphic disks (which are controlled by area bounds); then, the
motion in finite space is performed and the moving moduli
spaces are compactified as usual.
Finally, any compact motion is divided up into a
finite sequence of small motions, all of which lead to isomorphic objects.

\end{enumerate}

The preceding (1--5) ensure well-definedness of the Fukaya category.
We need more to apply the Fukaya-Oh theorem relating to the
Morse category.  The issue is that the theorem as stated only
applies to global graphs $\Gamma_{df}$ with $f:B\rightarrow \bR,$
whereas standard branes are graphs only over open sets $U.$
The main thing to prove is that the bounds ensuring compact
moduli spaces do not break down under the dilation from
$\Gamma_{d\vec{f}}$ to $\Gamma_{\epsilon d\vec{f}}.$

\subsection{Statement of results}

The arguments of this section are the ideas behind the theorem, stated
formally here:
\begin{theorem}\cite{NZ,N}
Let $B$ be a real analytic manifold.
There is a quasiequivalence of $A_\infty$ categories
$$
Sh_c(B)\cong TrFuk(T^*B)
$$
which sends $i_*\bC_U$ to $\Gamma_{d\log m},$ where $m$ is a defining
function for $\partial U.$
Taking cohomology $H^0$ leads to a quasi-equivalence of
bounded, derived categories
$$
DSh_c(B)\cong DFuk(T^*B).
$$
\end{theorem}
The theorem makes sense for noncompact $B$ if we restrict
to objects with compact support.
The theorem also makes sense if we fix a conical Lagrangian
inside $T^*B$ and consider the subcategory of constructible
sheaves with singular support contained in it.
Applied to the conical Lagrangian $\LS\subset T^*M_\bR$
defined in \eqref{eqn:LS-M}, 
we have
\begin{eqnarray}
DSh_{cc}(M_\bR;\LS) &\cong & D Fuk(T^*M_\bR;\LS).
\end{eqnarray}

\section{Coherent-Constructible Correspondence}
\label{sec:ccc}

The results of Sections \ref{sec:tduality} and \ref{sec:micro} produce a
functor between coherent and constructible sheaves that goes
through the Fukaya category.  In this section we cut out the middleman
and construct the coherent-constructible correspondence directly.
In so doing, we discover that we have stumbled upon a categorification of Morelli's
description of the equvariant K-theory of a toric variety \cite{M}, which we now describe.

We use the notation in Section \ref{subsec:statement-fltz}
and Section \ref{subsec:tobj}.
Let $X_\Si$ be a smooth projective toric variety
defined by a complete fan $\Si\subset N_\bR$.
Let $D_\vc$ be an ample $T$-divisor, so that
$\cO_{X_\Si}(D_\vc)$ is an equivariant ample
line bundle. Let $1_{\tri_{\vc}}:M_\bR  \to \bR$ denote the indicator
function on the moment polytope $\tri_{\vc}$
of $\cO_{X_\Si}(D_\vc)$. Then
$\cO_{X_\Si}(D_\vc)\mapsto 1_{\tri_{\vc}}$
defines a homomophism
\begin{equation}\label{eqn:IT}
I_T: K_T(X_\Si) \to L_M(M_\bR)
\end{equation}
between abelian groups, where $L_M(M_\bR)$ is the group
of functions on $M_\bR$ generated over $\bZ$
by the indicator functions of convex lattice
polyhedra.  Morelli \cite{M} has the following characterization
of the image of $I_T$. Define the group of polyhedra germs
$S_M(M_\bR)$ to be the abelian group generated
by rational convex cones in $M_\bR$, with relations
$$
[\sigma\cup \tau]=[\sigma]+ [\tau]-[\sigma\cap \tau]
$$
where $\sigma, \tau, \sigma\cup \tau$ are rational
convex cones. Then $S_M(M_\bR)$ is the group
of germs of functions in $L_M(M_\bR)$ at the origin (or at any
point in  $M$). Let $S_\Si(M_\bR)$ be the subgroup of $S_M(M_\bR)$
generated by $\{ \sigma^\vee\mid \sigma\in \Si\}$. Then
the image $L_\Si(M_\bR)$ of $I_T$ is the subgroup
of $L_M(M_\bR)$ consisting
of functions whose germ at any point $m\in M$
lies in $S_\Sigma(M_\bR)$, and $I_T$ defines
an isomorphism
\begin{equation}\label{eqn:KT}
K_T(X_\Si)\cong L_\Si(M_\bR)
\end{equation}
of abelian groups.

The coherent-constructible correspondence
\begin{equation}\label{eqn:cccT}
D Coh_T (X_\Si) \cong  DSh_{cc}(M_\bR;\LS)
\end{equation}
can be viewed as categorifcation of Morelli's results
\eqref{eqn:KT}.
%
In \cite{FLTZ}, the authors defined a functor
$\kappa: D Coh_T(X_\Si)\to Sh_{cc}(M_\bR;\LS)$
and proved that it is a quasi-equivalence
of dg categories. Taking the cohomology $H^0$
gives \eqref{eqn:cccT}. We outline the argument
in the remainder of this section.

Recall that $D Coh_T(X_\Si)$ is generated
by equivariant ample line bundles, so $\kappa$
is determined by its restriction to equivariant ample line
bundles. Given an equivariant ample line
bundle $\cO_{X_\Si}(D_\vc)$, which is
an object in $D Coh_T(X_\Si)$, let
$i:\tri^\circ_\vc\hookrightarrow M_\bR$ be the
inclusion map. Then $i_! \bC_{\tri^\circ_\vc}[n]$
is a constructible sheaf on $M_\bR$ with
compact support, where $n=\dim_\bC X_\Sigma$,
and $SS (i_! \bC_{\tri^\circ_\vc})\subset \LS$,
so it is an object of $Sh_{cc}(M_\bR;\LS)$.
We want to show that
$\cO_{X_\Si}(D_\vc)\mapsto i_! \bC_{\tri^\circ_\vc}[n]$
defines a functor
$\kappa: D Coh_T(X_\Si)\to Sh_{cc}(M_\bR;\LS)$
between dg categories.
To verify this, it is more convenient to express
$\kappa$ in terms of different generators (which do not
actually belong to $D Coh_T(X_\Si)$).
For any pair $(\chi,\sigma)\in M\times \Sigma$, we
will define a constructible sheaf $\Theta(\chi,\sigma)$
on $M_\bR$ and an equivariant quasi-coherenet sheaf
$\Theta'(\chi,\sigma)$ on $X_\Sigma$.

\begin{enumerate}
\item Given $(\chi,\sigma)\in M\times\Sigma$,  let
$(\chi+\sigma^\vee)^\circ$ be the interior of
the translated dual cone $\chi+ \sigma^\vee\subset M_\bR$.
Define
$$
\Theta(\chi,\sigma)= i_! \bC_{(\chi+\sigma^\vee)^\circ}.
$$
Then $\Theta(\chi,\sigma)$ is a constructible
sheaf on $M_\bR$ whose microlocal support
is contained in $\LS$, but its support is not compact.
Therefore $\Theta(\chi,\sigma)$ is an object
of $Sh_c(M_\bR;\LS)$.
Let $\langle \Theta \rangle$ denote the full triangulated
dg subcategory of $Sh_c(M_\bR;\LS)$ generated by
$$
\{ \Theta(\chi,\sigma)\mid (\chi,\sigma)\in M\times \Sigma\}.
$$
The hom space between any two of these generators is simple:
the computation can be reduced to computing
relative cohomology groups
of (contractible set, point) or
(contractibe set, empty set),
which is either $0$, or $\bC$ (at degree zero).

\item Given $\sigma\in\Sigma$, $X_\sigma=\Spec\bC[\sigma^\vee\cap M]$
is an affine open subvariety of $X$. Let
$\cO_\sigma(\chi)$ be the equivariant quasicoherent
sheaf on $X_\sigma$ corresponding to the $M$-graded
$\bC[\sigma^\vee\cap M]$-module freely generated
by $\chi\in M$. Define
$$
\Theta'(\chi,\sigma)= j_{\sigma *}\cO_\sigma(\chi)
$$
where $j_{\sigma *}:X_\sigma \hookrightarrow X$ is
the open embedding. Then $\Theta'(\chi,\sigma)$
is an equivariant quasi-coherent sheaf on $X_\Si$.
Let $\langle \Theta'\rangle$ denote the full triangulated
dg subcategory of $\cQ_T(X_\Si)$ (localization of the dg category of equivariant
quasi-coherent sheaves on $X_\Si$ with respect
to acyclic complexes) generated by
$$
\{ \Theta'(\chi,\sigma)\mid (\chi,\sigma)\in M\times \Sigma\}.
$$
The hom space between any two of these generators is
simple: the computation can be reduced to computing
$\Ext^*$ between two equivariant line bundles
on the affine toric variety $X_\sigma$, which is
either $0$,  or $\bC$ (at degree zero).
By \v{C}ech resolution we see that
any equivariant line bundle is in the
subcategory $\langle\Theta'\rangle$, so
$$
D Coh_T(X_\Si)\subset \langle \Theta'\rangle\subset \cQ_T(X_\Si).
$$

\item Comparing the calculations in (1) and (2),
we conclude that
$\Theta'(\chi,\sigma)\mapsto \Theta(\chi,\sigma)[n]$
defines a quasi-equivalence
$\langle \Theta' \rangle \cong \langle \Theta\rangle$.
Moreover, under this quasi-equivalence an ample line bundle
$\cO_{X_\Si}(D_\vc)$ is mapped to
$i_! \bC_{\tri_{\vc}^\circ}[n]$.
Therefore we have a full embedding
$$
\kappa:D Coh_T(X_\Si)\to Sh_{cc}(M_\bR;\LS).
$$

\item In \cite{FLTZ}, it is shown that $\kappa$ is essentially surjective, i.e. that every sheaf
$F \in Sh_{cc}(M_\bR;\LS)$ is quasi-isomorphic to one of the form $\kappa(\cG)$.  This requires a careful induction on the ``height'' of $F$, which is a measure of the complexity of its singular support.

\item The argument can be extended to any complete toric variety, including singular and non-projective varieties.  The category $D Coh_T(X_\Si)$ is replaced by $\Perf_T(X_\Si)$, the dg category of ``perfect complexes,'' which are by definition bounded complexes of equivariant vector bundles.

\item  The techniques of \cite{FLTZ} produce a full embedding
$$
\bar{\kappa}: DCoh(X_\Sigma) \to Sh_c(T_\bR^\vee,\LS/M)
$$
which makes the following square commute:
$$
\begin{CD}
D Coh_T(X_\Si) @>{\kappa}>> Sh_{cc}(M_\bR;\LS)\\
@V{f}VV  @V{p_!}VV\\
D Coh(X_\Si) @>{\bar{\kappa}}>> Sh_{cc}(T_\bR^\vee;\LS/M).
\end{CD}
$$
where $f$ is obtained by forgetting
the $T$-equivariant structure, and $p_!$
is induced by the natural projection
$p: M_\bR \to T^\vee_\bR=M_\bR/M$.
Presumably $\bar{\kappa}$ is always an equivalence,
though as of this writing we do not have a proof of this
fact. (For example, $\bar{\kappa}$ is an equivalence
in the examples in Section \ref{subsec:Pone} and
Section \ref{subsec:fano-surface} below.)
\end{enumerate}

\section{Examples} \label{sec:exs}

\subsection{Taking the mapping cone} \label{subsec:Pone}

On $B=S^1=\bR/\bZ$, let $p=0$ and $U=(0,1)$. We have $S^1=\{p\}\cup U$. There is an exact triangle in the derived category
$$\bC_{S^1} \to i_{U*}\bC_{U}\to \cO_p\to $$
where $i_U: U\hookrightarrow S^1$ is the embedding map and $\cO_p$ is the skyscraper sheaf at $p$.

Let $F_p$ be the fiber of $T^*S^1\to S^1$ at $p$, $L_{U*}$ be the standard Lagrangian over the interval $U$. The microlocalization functor $$\mu: Sh_c(S^1)\to TrFuk(T^*S^1)$$ takes the above exact cone to
$$B\to L_{U*} \to F_p\to.$$
The brane $F_p$ and $L_{U*}$ are equipped with the trivial gradings, and $B$ is simply the Lagrangian brane supported at the zero section $S^1$ (with the trivial grading).

Under this functor, $K$-theory classes
$K(i_{U*}\bC_{U})=K(\bC_{S^1}\oplus \cO_p)$ map to $K(L_{U*})=K(B\oplus F_p)$, which reflects the dilation limit
$$\lim_{\epsilon\to 0} \epsilon L_{U*}= F_p \cup S^1.$$

Taking the mapping cone of
the morphism $B\to L_{U*}$ is geometrically interpreted as taking symplectic surgery,
as depicted in the Figure 5.\footnote{``Cone = Surgery'' is the mantra found from the
lessons of \cite{FOOO,S3,ST}.}
\begin{figure}[h]
\begin{center}
\includegraphics[scale=0.4]{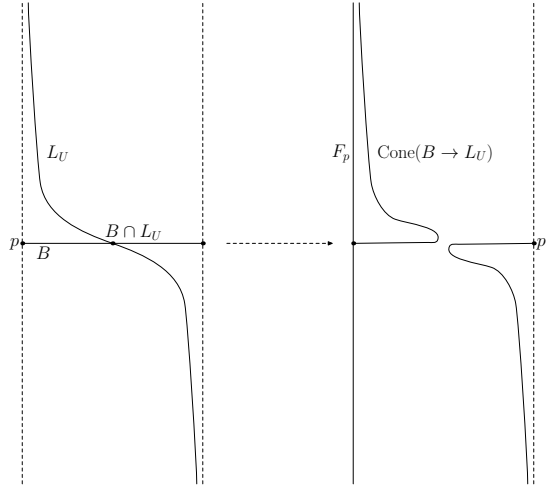}
\end{center}
\caption{The mapping cone over the morphism
$B\to L_{U*}$ in the Fukaya category is obtained by taking the symplectic surgery at the intersection point (left). The resulting brane ${\tt Cone} (B\to L_{U*})$ (right) is isotopic to the fiber brane $F_p$. Notice in both pictures the cotangent bundle $T^* S^1$ is cut along $F_p$ for illustration.}
\end{figure}

\begin{remark}
When $X_\Si=\bP^1$, $M_\bR\cong S^1$ and $\LS/M=F_p\cup B$.
We have the (non-equivariant) mirror symmetry
$$
DFuk(T^*(M_\bR/M);\LS/M)\cong D Coh(\bP^1).
$$Let $D_0$ be one of the $T$-invariant divisors. The
above exact triangle comes from
$$\cO_{D_0}[-1] \to \cO(-1)\to \cO \to.$$
\end{remark}

\subsection{Toric Fano surfaces} \label{subsec:fano-surface}
 Let $X_\Sigma$ be a smooth projective toric Fano surface,
so that $\Sigma$ is one of the five fans
in the first row of Figure 3. We have
$$
DCoh_T(X_\Sigma)\cong DSh_{cc}(M_\bR;\LS),\quad
DCoh(X_\Sigma)\cong DSh_c(T_\bR^\vee;\LS/M).
$$
Moreover, in this case
there is a constructible (but possibly
non-Whitney) stratification of $T_\bR^\vee$
(see Figure 6), such that each stratum
is contractible, and any object in $Sh_c(T_\bR^\vee;\LS/M)$
is constant on each stratum.
This stratification coincides with
 the one defined by Bondal \cite{Bondal}.

\begin{figure}[h]
\begin{center}
\psfrag{F0}{$\bF_0=\bP^1\times \bP^1$}
\psfrag{P2}{$\bP^2$}
\psfrag{B1}{$B_1=\bF_1$}
\psfrag{B2}{$B_2$}
\psfrag{B3}{$B_3$}
\includegraphics[scale=0.6]{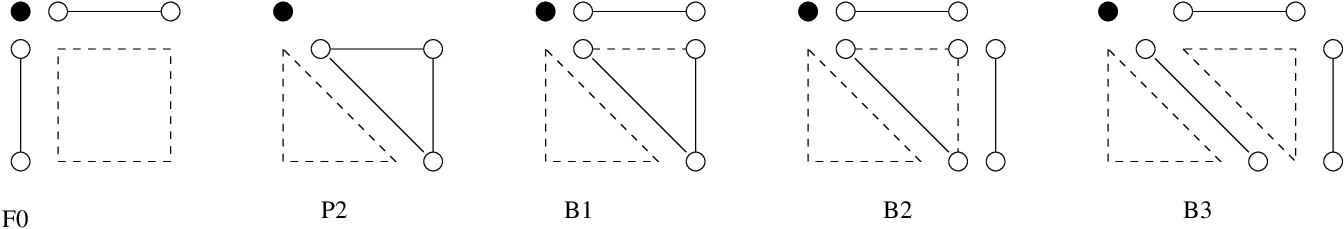}
\end{center}
\caption{$Sh_{c,\LS}(T_\bR^\vee)$ for toric Fano surfaces}
\end{figure}

\subsection{Hirzebruch surfaces}
\label{subsec:hirz}
The Hirzebruch surface $\bF_m$ ($m\geq 0$) is defined by a fan
$\Sigma_m$ spanned by four $1$-cones $\rho_i=\bR_{\geq 0}v_i$, where
$$
v_1=(1,0),\ v_2=(0,1),\ v_3=(-1,-m),\ v_4=(0,-1).
$$
See Figure 3 for the fans $\Sigma_m$, $m=0,1,2,3$.
Let $D_i$ be the $T$-divisor corresponding to $v_i$. Then
$D_2, D_4$ are sections of the projection $\bF_m\to \bP^1$, and
$D_1, D_3$ are fibers. We have $D_2\cdot D_2 =m = -D_4\cdot D_4$.
The fan data determines a conical Lagrangian $\Lambda_m:=\Lambda_{\Sigma_m}$
in $T^*M_\bR$. See Figure 4 for the conincal Lagrangians
$\Lambda_m$, $m=0,1,2,3$.

Let us discuss the various categories for this family of examples.

\subsubsection{Category of equivariant coherent sheaves $Coh_T(\bF_m)$}
This category is  generated by line bundles
$$
\cO_{\bF_m}(D_\vc),
$$
where $D_\vc=\sum_{i=1}^4 c_i D_i$.
This category is also generated by \emph{ample} line bundles.

\subsubsection{ Equivariant category of constructible sheaves $Sh_{cc}(M_\bR;\Lambda_m)$}
This category consists of compactly-supported constructible sheaves whose singular supports
lie in $\Lambda_m$. It is generated by costandard sheaves
$i_{\tri^\circ_\vc\ !} \bC_{\tri^\circ_\vc}[2]$ for
all ample divisors $D_\vc$. When $D_\vc$ is ample, i.e.
$$
c_2+c_4>0,\ c_1 +c_3 > m c_4,
$$
the polytope $\tri_\vc$ is a polytope on the $M_\bR$ plane with four vertices
$$
(-c_1,-c_2),\ (c_3 + mc_2, -c_2),\ (c_3- m c_4, c_4),\ (-c_1, c_4).
$$
The map $i_{\tri^\circ_\vc}$ is the inclusion $\tri^\circ_\vc \hookrightarrow M_\bR$
of the interior.

In particular, the anti-canoncial divisor of $\bF_m$ is given by
$$
-K_{\bF_m}= D_1+ D_2+ D_3 + D_4.
$$
It is ample when $m=0,1$, and is numerically effective when $m=2$.
It is not numerically effective when $m>2$.

\begin{figure}[h]
\begin{center}
\psfrag{(0,0)}{\tiny $(0,0)$}
\psfrag{(0,1)}{\tiny $(0,1)$}
\psfrag{(1,1)}{\tiny $(1,1)$}
\psfrag{(2,0)}{\tiny $(2,0)$}
\includegraphics[scale=0.6]{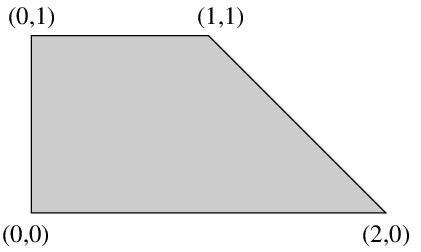}
\end{center}
\caption{The polytope $\tri_{(0,0,2,1)}$ for $\bF_1$.
The constructible sheaf corresponding to the ample bundle $\cO_{\bF_1}(2D_3+D_4)$ is the costandard constructible sheaf over this polytope (i.e. $i_{\tri^\circ_{(0,0,2,1)}\ !} \bC_{\tri^\circ_{(0,0,2,1)}} [2])$.}
\end{figure}

\subsubsection{Equivariant Fukaya category $Fuk(T^*M_\bR;\Lambda_m)$ associated to the universal cover of the mirror
$\tilde Y^\vee \cong T^*M_\bR$  (symplectically)} The category $Fuk(T^*M_\bR;\Lambda_m)$ consists of
Lagrangian branes whose conical limits are subsets of $\Lambda_m$.
It is generated by costandard Lagrangian branes $\bL_{h, \vc}$ over $\tri_\vc$
for all ample divisors $D_\vc$.

\subsubsection{Tri-Equivalence}
The three categories above are equivalent, by \cite{FLTZ}:
$$
Coh_T(\bF_m)\stackrel{\cong}{\longrightarrow}
Sh_{cc}(M_\bR;\Lambda_m)\stackrel{\cong}{\longrightarrow} Fuk(T^*M_\bR;\Lambda_m).
$$
The correspondence between a generating class of objects is
$$
\cO_{\bF_m}(D_\vc)\mapsto \bL_{\vc, h}\mapsto i_{\tri^\circ_\vc\ !} \bC_{\tri^\circ_\vc}[2].
$$

%

\subsubsection{Relative Fukaya category $Fuk((\bC^*)^2, W_m^{-1}(0))$}
Let $\vt=(t_1, t_2,t_3,t_4)$ be chosen such that
$D_\vt=\sum_{i=1}^4 t_i D_i$ is an ample divisor. Then
the Poincar\'{e} dual of this divisor class is a K\"{a}hler class
in $H^2(\bF_m;\bR)$.
The superpotential $W_m: (\bC^*)^2\rightarrow \bC$ of the
Hori-Vafa mirror of $\bF_m$ is given by
\begin{equation}\label{eqn:Wm}
W_m=e^{- t_1} z_1+e^{-t_2} z_2+ e^{-t_3} z_1^{-1}z_2^{-m}+e^{-t_4}z_2^{-1}.
\end{equation}

An object in $Fuk((\bC^*)^2,W_m^{-1}(0))$ is an exact compact
Lagrangian submanifold with boundary in the complex
hypersurface $W_m^{-1}(0)$.
Abouzaid defined a subcategory
$\cT Fuk ((\bC^*)^2,W_m^{-1}(0))$ of tropical Lagrangian
sections of $Fuk((\bC^*)^2,W_m^{-1}(0))$, and showed that
$\cT Fuk((\bC^*)^2,W_m^{-1}(0))$ is equivalent to the category
of line bundles on $\bF_m$. Since line bundles generate
$DCoh(X)$, he concluded that
\begin{equation}\label{eqn:TFuk}
DCoh(X) \cong D^\pi \cT Fuk((\bC^*)^2, W_m^{-1}(0))
\stackrel{i_{\cT}}{\hookrightarrow} D^\pi Fuk((\bC^*)^2, W_m^{-1}(0))
\end{equation}
where $i_{\cT}$ is embedding of a full subcategory.
(Indeed, Abouzaid proved a statement like \eqref{eqn:TFuk}
for any smooth projective toric variety.)

To describe a tropical Lagrangian
section, we consider the tropical limit $M_\infty$ of $W_m^{-1}(0)$
defined by the ample divisor $D_\vt$
\cite{Mi}. The image
of $M_\infty$ under the logarithm map
$\log: (\bC^*)^2\to \bR^2$ is a graph known as the tropical amoeba
(see Figure 8).
\begin{figure}[h]
\begin{center}
\psfrag{F0}{$\bF_0=\bP^1\times \bP^1$}
\psfrag{P2}{$\bP^2$}
\psfrag{B1}{$B_1=\bF_1$}
\psfrag{B2}{$B_2$}
\psfrag{B3}{$B_3$}
\psfrag{F2}{$\bF_2$}
\psfrag{F3}{$\bF_3$}
\includegraphics[scale=0.6]{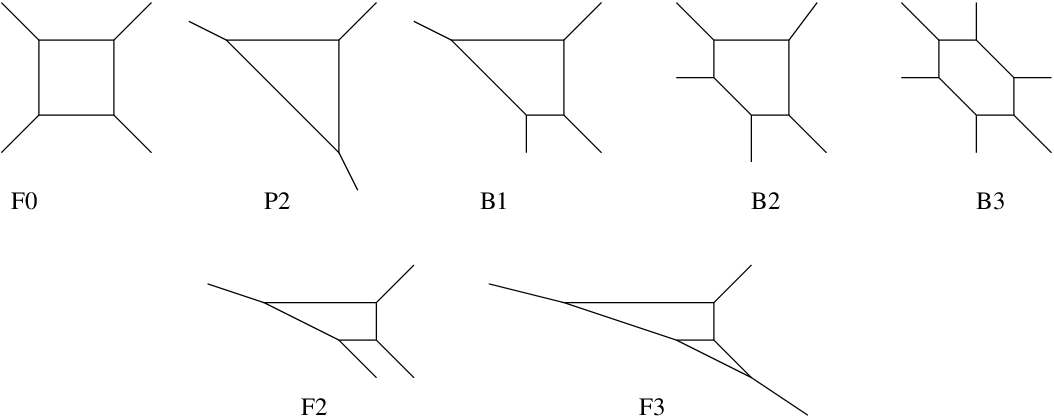}
\end{center}
\caption{Tropical amoeba of some toric surfaces}
\end{figure}
The tropical amoeba divides
$\bR^2$ into connected components.
When $m=0,1,2$, there is a unique bounded component
which can be identified with the moment polytope
of the anti-ample divisor $D_{-\vt}$. When $m> 2$,
there are more than one bounded components, one
of which can be identified with the moment
polytope of $D_{-\vt}$.
In the first row of Figure 8, we take $D_{\vt}$ to be
the anti-canonical divisor, which is ample; in
the second row, we take $D_{\vt}= D_1+ D_2$.
An object in $\cT Fuk((\bC^*)^2,M_\infty)$ is
a section of $\log^{-1}(\tri_{-\vt})\to \tri_{-\vt}$
with boundary in $M_\infty$.  Let $e_i$
be the face of $\tri_{-\vt}$ which corresponds
to $D_i$. Then $e_i$ is a closed interval parallel to the line $\rho_i^\perp\subset M_\bR$.
Let $e_i^\circ$ be the interior of $e_i$. Then
$$
\log^{-1}(e_i^\circ) \cap M_\infty = T e_i^{\circ}/ (M\cap \rho_i^\perp).
$$

Given an object
$$
\bL_{\vc,h}=\{ (\Phi_h\circ j_0(y),y)\mid y\in N_\bR \} \subset M_\bR \times N_\bR
$$
in $Fuk(T^*M_\bR;\Lambda_m)$, we will associate
a tropical Lagrangian section. Define
$$
L_{\vt,\vc,h}=\{ (\Phi_h\circ j_0(y), -\Phi_\vt\circ j_0(y))\mid y\in N_\bR\}\subset M_\bR\times M_\bR
$$
where $\Phi_\vt:\bF_m\to \bR^2$ is a moment map of the ample divisor $D_\vt$.
Then $L_{\vt,\vc,h}$ is a section of $M_\bR \times \tri_{-\vt}^\circ \to \tri_{-\vt}^\circ\,$;
its closure is
$$
\overline{L}_{\vt,\vc,h} =\{ (\Phi_h(x), -\Phi_\vt(x))\mid x\in \bF_m \}
$$
which is a section of $M_\bR\times \tri_{-\vt}\to \tri_{-\vt}$.
Moreover, the boundary  of $\overline{L}_{\vt,\vc,h}$ is contained in
$$
\bigcup_{i=1}^4 (M  + \rho_i^\perp)\times F_i,
$$
so $\overline{L}_{\vt,\vc,h}$ descends to a section $\cT_{\vt,\vc,h}$
of $\log^{-1}(\tri_{-\vt})\to \tri_{-\vt}$ with boundary in $M_\infty$.

Abouzaid showed that one can smooth $\cT_{\vt,\vc,h}$ to
obtain a tropical section $L'_{\vc,h}$ which is compact Lagrangian with
boundary in $W_m^{-1}(0)$, a smooth hypersurface
very close to $M_\infty$ (up to scaling). Then $L'_{\vc,h}$
is an object in $\cT Fuk((\bC^*)^2,W_m^{-1}(0))$.

It is expected that the embedding $i_{\cT}$ in \eqref{eqn:TFuk}
is a  quasi-equivalence when $m=0,1$. (More generally,
one expects such a quasi-equivalence for any smooth projective
toric {\em Fano} variety. The case for $m=2$ is also a quasi-equivalence here.)
When $m>2$, there are other
bounded regions giving rise to Lagrangian
sections which do not correspond to objects
in $D Coh(\bF_m)$, so $i_{\cT}$ cannot be
a quasi-equivalence.

\subsubsection{Fukaya-Seidel category $FS((\bC^*)^2, W_m)$}
Let the superpotential $W_m$ be defined by \eqref{eqn:Wm} as above.
When $m=0,1$, $\bF_m$ is Fano, we may take $D_\vt$
to be the anti-canonical class,
i.e., $t_1=t_2=t_3=t_4=1$. So
$$
W_1= e^{-1}(z_1+ z_2+ z_1^{-1} z_2^{-1}+ z_2^{-1}).
$$
There are four critical values of $W_1$, namely,
$\lambda_1>0,\lambda_3<0$ on the real line,
and two mutually conjugate $\lambda_0$ and $\lambda_2$ with $\mathrm{Im} \lambda_0>0$.
To each critical value, there is a unique critical point in $(\bC^*)^2$. Choose a generic fiber
$W_1^{-1}(\lambda)$ and choose four paths from $\lambda$ to the critical values.
The vanishing cycles associated to these paths are the objects of
the Fukaya-Seidel category $FS((\bC^*)^2, W_1)$.
For example, Figure 9 describes a particular layout of paths and thus defines a Fukaya-Seidel category.
Denote the vanishing cycle associated to the critical point $\lambda_i$ by $L_i$.
\begin{figure}[h]
\begin{center}
\psfrag{lambda1}{\footnotesize $\lambda_1$}
\psfrag{lambda2}{\footnotesize $\lambda_2$}
\psfrag{lambda3}{\footnotesize $\lambda_3$}
\psfrag{lambda}{\footnotesize $\lambda$}
\psfrag{lambda0}{\footnotesize $\lambda_0$}
\includegraphics[scale=0.3]{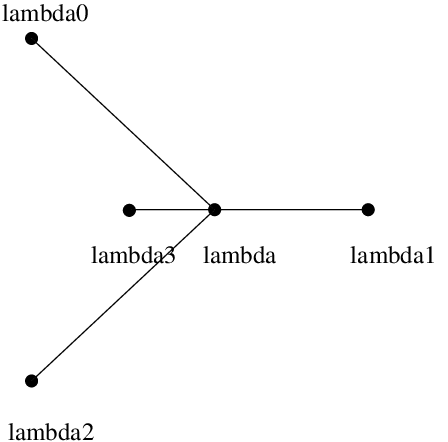}
\end{center}
\caption{A choice of paths defining the Fukaya-Seidel category for
$((\bC^*)^2,W_1)$. The category consists of four vanishing cycles in
$W_1^{-1}(\lambda)$. A different homotopy type of these paths gives
rise to the mutation of objects, as described in \cite{S2},
\cite{S2.5}.}
\end{figure}

 Auroux, Katzarkov and Orlov \cite{AKO} prove the quasi-equivalence
$DCoh(\bF_1)\cong DFS((\bC^*)^2, W_1)$.
Let $\pi$ be the blow-up map $\bF_1\to \bP^2$ and let $E=D_4$ be the exceptional curve in $\bF_1$.
Object by object, they show the mirror correspondence
$$
L_0\sim \cO_{\bF_1},\ L_1\sim  \pi^*(T_{\bP^2}(-1)),\ L_2\sim \pi^*(\cO_{\bP^2}(1)),\ L_3\sim \cO_E.
$$
The full strong exceptional collection $\{\cO_{\bF_1},\ \pi^*(T_{\bP^2}(-1)),\ \pi^*(\cO_{\bP^2}(1)),\ \cO_E\}$
implies $DFS((\bC^*)^2, W_1)\cong DCoh(\bF_1)$.

In general, we may take $D_\vt= D_1+D_2$, i.e.,
$t_1=t_2=1$, $t_3=t_4=0$. Then
$$
W_m = e^{-1} z_1 + e^{-1} z_2 + z_1^{-1} z_2^{-m}+ z_2^{-1}.
$$
There are $m+2$ critical points $\lambda_0,\lambda_1,\ldots, \lambda_{m+1}$
of $W_m$. Then $DFS((\bC^*)^2,W_m)$ is generated
by $m+2$ vanishing cycles $L_0, L_1,\ldots, L_{m+1}$. Auroux, Katzarkov
and Orlov showed that, in the limit
$$
W_{m,t_3} =e^{-1} z_1 + e^{-1} z_2 + e^{-t_3} z_1^{-1} z_2^{-m} + z_2^{-1},\quad t_3\to \infty,
$$
four of the critical points stay in a bounded region
(we may assume they are $\lambda_0,\ldots,\lambda_3$),
while the other $m-2$ critical points go to infinity.
They showed that
\begin{equation}\label{eqn:four}
DCoh(\bF_m) \cong D \cB FS((\bC^*)^2, W_m)
\stackrel{i}{\hookrightarrow} D FS((\bC^*)^2,W_m)
\end{equation}
where $\cB FS((\bC^*)^2, W_m)\subset FS((\bC^*)^2, W_m)$ is
the subcategory generated by $L_0, L_1, L_2, L_3$.
The inclusion $i$ in \eqref{eqn:four}  is not an equivalence when $m>2$.

We may take $\lambda=0$. Then
the tropical Lagrangian section $L'_{\vc,h}$ is
a circle in $W_m^{-1}(0)$. In particular, the
zero section can be identified with $L_0$.
$L'_{\vc,h}$ is contained in a bounded region.
It is expected that $L'_{\vc,h}$
is contained in the subcategory
$D \cB FS((\bC^*)^2, W_m)$,
and that the Lagrangian sections
$L'_{\vc,h}$'s also generate
$D\cB FS((\bC^*)^2, W_m)$.

\begin{remark}
Here the superpotential $W_m$ corresponds to
the leading order of of the full potential in Fukaya-Oh-Ohta-Ono \cite{FOOO1}
and Auroux \cite{Au2}.

In \cite{FOOO1}, the full potential
\begin{equation}\label{eqn:PO}
\frak{PO} =\frak{PO}_0 +\textup{ higher order terms}
\end{equation}
counts holomorphic disks of Maslov index 2 in $\bF_m$ with
boundary in a Lagrangian torus fiber of the moment map.
The higher order terms in \eqref{eqn:PO} come from disks
with sphere bubbles; such index 2 disks do not
exist in the Fano case ($m=0,1$), but
exist and contribute to $\frak{PO}$ when $m\geq 3$.
Fukaya-Oh-Ohta-Ono showed that $\frak{PO}$
has four critical points in $Y^\vee\subset T^\vee
=\Spec\bC[z_1, z_1^{-1}, z_2, z_2^{-1}]$ (or rather
$\Spec (\Lambda_{nov}\otimes \bC)[z_1, z_1^{-1}, z_2, z_2^{-1}]$)
\cite[Example 7.2]{FOOO1}.
In a forthcoming sequel of \cite{FOOO1, FOOO2},
Fukaya-Oh-Ohta-Ono prove that the number
of critical points of $\frak{PO}$  \footnote{When there are degenerate
critical points, one counts with multiplicities.} in $Y^\vee$  is equal to
$\dim_{\bQ} H^*(X,\bQ)$  for any smooth projective toric variety $X$.

\end{remark}



\begin{thebibliography}{CK}

\bibitem{Ab1} M. Abouzaid,
``Homogeneous coordinate rings and mirror symmetry for toric varieties,''
Geom. Topol. {\bf 10} (2006), 1097--1157.

\bibitem{Ab2} M. Abouzaid,
``Morse Homology, Tropical Geometry, and
Homological Mirror Symmetry for Toric Varieties'',
{\tt arXiv:math/0610004}.

\bibitem{AP} D. Arinkin and A. Polishchuk,
``Fukaya category and Fourier transform,''
in {\em Winter School on Mirror Symmetry, Vector Bundles and Lagrangian Submanifolds,}
AMS/IP Stud. Adv. Math. {\bf 23} (2001) 261--274.

\bibitem{AB} M. F. Atiyah and R. Bott, ``The moment map and equivariant cohomology,''  Topology {\bf 23} (1984), no. 1, 1--28.

\bibitem{Au1} D. Auroux, ``Mirror symmetry and T-duality in the complement of
the anticanonical divisor,"  {\tt arXiv:0706.3207}.

\bibitem{Au2} D. Auroux, ``Special lagrangian fibrations, wall-crossing, and mirror
symmetry," {\tt arXiv:0902.1595}.

\bibitem{AKO} D. Auroux, L. Katsarkov, D. Orlov,
``Mirror symmetry for weighted projective planes and their
noncommutative deformations'', Ann. of Math. (2)  {\bf 167}  (2008),
no. 3, 867--943.

\bibitem{AKO2} D. Auroux, L. Katsarkov, D. Orlov, ``Mirror symmetry for del Pezzo surfaces:
vanishing cycles and coherent sheaves'',  Invent. Math.  {\bf 166}
(2006), no. 3, 537--582.

\bibitem{Bondal} A. Bondal, ``Derived categories of toric varieties,'' in
{\em Convex and Algebraic geometry, Oberwolfach conference reports,
EMS Publishing House} {\bf 3} (2006) 284--286.

\bibitem{Ch1} C.-H. Cho,
``Products of Floer cohomology of torus fibers in toric Fano manifolds,''
Comm. Math. Phys. {\bf 260} (2005), no. 3, 613--640.

\bibitem{Ch2} C.-H. Cho,
``On the counting of holomorphic discs in toric Fano manifolds,''
{\tt arXiv:math/0604502}.

\bibitem{CO} C.-H. Cho, Y.-G. Oh,
``Floer cohomology and disc instantons of Lagrangian torus fibers
in Fano toric manifolds'',
Asian J. Math. {\bf 10} (2006), no. 4, 773--814.

\bibitem{CJS} R.L. Cohen, J.D. Jones, G.B. Segal,
``Floer's infinite-dimensional Morse theory and homotopy theory,''
{\em The Floer memorial volume}, 297--325,
Progr. Math., {\bf 133}, Birkh\"{a}user, Basel, 1995.

\bibitem{CK} D. Cox and S. Katz, {\em Mirror symmetry and algebraic geometry},
Mathematical Surveys and Monographs {\bf 68}, AMS, 1999.

\bibitem{CL1} K. Chan and N.C. Leung,
``Mirror symmetry for toric Fano manifolds via SYZ transformations'',
{\tt arXiv:0801.2830}.

\bibitem{CL2} K. Chan and N.C. Leung,
``On SYZ mirror transformation'',
{\tt arXiv:0808.1551}.

\bibitem{E} D. Eisenbud,
``Homological algebra on a complete intersection, with an application to group representations,''
Trans. Am.  Math. Soc. {\bf 260} (1980) 35--64.


\bibitem{Fa} B. Fang,
``Homological mirror symmetry is T-duality for $\mathbb P^n$,''
{\tt arXiv:0804.0646}.

\bibitem{FLTZ} B. Fang, C.-C.M. Liu, D. Treumann, E. Zaslow,
``
T-Duality and Equivariant Homological Mirror Symmetry for Toric Varieties,''
{\tt arXiv:0811.1228v1}.

\bibitem{Fl} A. Floer, ``Witten's complex and infinite-dimensional Morse theory,''
J. Differential Geom. {\bf 30} (1989), no. 1, 207--221.

\bibitem{F} K. Fukaya ``Morse Homotopy and its Quantization,''
AMS/IP Studies in Advanced Mathematics {\bf 2} (1997) 409--440.

\bibitem{FO} Fukaya and Y.-G. Oh, ``Zero-loop open strings
in the cotangent bundle and Morse homotopy,'' Asian J. Math {\bf 1}
(1997) 96--180.

\bibitem{FOOO} K. Fukaya, Y.-G. OH, H. Ohta, K. Ono, ``Lagrangian
Floer Theory -- Anomaly and Obstruction,'' Kyoto preprint Math 00-17, 2000.

\bibitem{FOOO1}K. Fukaya, Y.-G. OH, H. Ohta, K. Ono,
``Lagrangian Floer theory on compact toric manifolds I'',
{\tt arXiv:0802.1703}.

\bibitem{FOOO2}K. Fukaya, Y.-G. OH, H. Ohta, K. Ono
``Lagrangian Floer theory on compact toric manifolds II : Bulk deformations'',
{\tt arXiv:0810.5654}.

\bibitem{FSS}K. Fukaya, P. Seidel, and I. Smith,
``Exact Lagrangian submanifolds in simply-connected cotangent bundles,"
{\tt arXiv:math/07-1783}.

\bibitem{Fu}W. Fulton, {\em Introduction to toric varieties,}
Annals of Mathematics Studies {\bf 131},
Princeton University Press, 1993.

\bibitem{Gu} V. Guillemin,
``Kaehler structures on toric varieties,''
J. Differential Geom. {\bf 40}  (1994),  no. 2, 285--309.


\bibitem{HV} K. Hori and C. Vafa, ``Mirror Symmetry,'' {\tt hep-th/0002222.}

\bibitem{HIV} K. Hori, A. Iqbal, C. Vafa,
``D-branes and mirror symmetry,'' {\tt hep-th/0005247.}

\bibitem{I} V. Iacovino, ``A simple proof of a theorem of Fukaya and Oh,''
{\tt arXiv:0812.0129}.

\bibitem{K} B. Keller, ``Introduction to $A$-infinity algebras and
modules,''  Homology Homotopy Appl. {\bf 3}  (2001) 1--35.

\bibitem{Ko} M. Kontsevich,
``Homological algebra of mirror symmetry,''
Proceedings of the International Congress of Mathematicians, Vol. 1, 2
(Z\"{u}rich, 1994),  120--139, Birkh\"{a}user, Basel, 1995.

\bibitem{Ko2} M. Kontsevich, course at ENS, 1998,
{\tt
http://www.math.uchicago.edu/$\sim$mitya/langlands/kontsevich.ps}

\bibitem{KaSc} M. Kashiwara and P. Schapira, {\em Sheaves on Manifolds,}
Grundlehren der Mathematischen Wissenschaften {\bf 292}, Springer-Verlag
(1994).

\bibitem{KoSo} M. Konstevich and Y. Soibelman,  ``Homological mirror symmetry and torus fibrations,''
in {\em Symplectic geometry and mirror symmetry (Seoul, 2000)} World Sci. Publ. (2001) 203--263.

\bibitem{KT} Y. Karshon and S. Tolman,
``The moment map and line bundles over pre-symplectic toric
mainifolds,"  J. Differential Geom.  {\bf 38}  (1993),  no. 3, 465--484.

\bibitem{LYZ} N.-C. Leung, S.-T. Yau and E. Zaslow,
``From special Lagrangian to Hermitian-Yang-Mills via Fourier-Mukai transform,"
Adv. Theor. Math. Phys. {\bf 4} (2000) 1319--1341.

\bibitem{Mi} G. Mikhalkin,
``Decomposition into pairs-of-pants for complex algebraic surfaces,''
Topology {\bf 43} (2004), no. 5,  1035--1065.

\bibitem {M} R. Morelli, ``The K theory of a toric variety'',
Adv. Math.  {\bf 100}  (1993),  no. 2, 154--182.


\bibitem{N} D. Nadler,
``Microlocal branes are constructible sheaves'',
{\tt arXiv:math/0612399}.

\bibitem{NZ} D. Nadler, E. Zaslow,
``Constructible sheaves and the Fukaya category'',
J. Amer. Math. So. {\bf 22} (2009), 233-286.


\bibitem{Or} D. Orlov,
``Triangulated categories of singularities and D-branes in Landau-Ginzburg
models'', {\tt arXiv:math/0302304}.

\bibitem{SV} W. Schmid and K. Vilonen, ``Characteristic Cycles
of constructible sheaves,'' Invent. Math. {\bf 124} (1996) 451--502.

\bibitem{S} P. Seidel, {\em Fukaya Categories and Picard-Lefschetz Theory,}
Zurich Lectures in Advanced Mathematics. European Mathematical Society (EMS),
Z\"{u}rich, 2008.

\bibitem{S2} P. Seidel, ``Vanishing cycles and mutation",
Proc. 3rd European Congress of Mathematics (Barcelona, 2000), vol.
II. Prog. Math. {\bf 202}, 65--85.

\bibitem{S2.5} P. Seidel, ``More about vanishing cycles and mutation",
Symplectic Geometry and Mirror Symmetry, Proc. 4th KIAS
International Conference (Seoul, 2000), 429-465.

\bibitem{S3} P. Seidel, ``A long exact sequence for symplectic Floer
cohomology", Topology {\bf 42} (2003) 1003--1063.


\bibitem{ST} P. Seidel\ and\ R. Thomas,
``Braid group actions on derived categories of coherent sheaves,"
Duke Math. J. {\bf 108} (2001) 37--108.

\bibitem{Si} J.-C. Sikorav,
``Some properties of holomorphic curves in almost complex manifolds,''
{\em Holomorphic curves in symplectic geometry,} 165--189,
Progr. Math., {\bf 117}, Birkh\"{a}user, Basel, 1994.


\bibitem {SYZ} A. Strominger, S.-T. Yau, and E. Zaslow,
``Mirror symmetry is $T$-duality," Nuclear Phys. {\bf B479} (1996) 243--259.

\bibitem{U} K. Ueda, ``Homological mirror symmetry for toric del Pezzo surfaces,"  Comm. Math. Phys.  {\bf 264}  (2006),  no. 1, 71--85.

\bibitem{UY} K. Ueda and M. Yamazaki, ``Homological mirror symmetry for toric orbifolds of toric del Pezzo
surfaces", {\tt arXiv:math/0703267}.

\bibitem{Va} C. Vafa,
``Topological Mirrors and Quantum Rings'',
Essays on Mirror Manifolds.

\end{thebibliography}
\end{document}